\documentclass[aos,MSNbibl,dvips]{arximspdf}
\usepackage{graphicx}

\doi{10.1214/12-AOS1044} %kopijuoti is PTS
\volume{41}
\issue{2}
\pubyear{2013}
\firstpage{508}
\lastpage{535}

\makeatletter
\newcommand{\rrvert}{\vert}
\newcommand{\llvert}{\vert}
\newcommand{\eqref}[1]{(\ref{#1})}
\newcommand{\indep}{\perp\hspace*{-6.2pt}\perp}
\newcommand{\indexset}{\mathbb{I}}
\newtheorem{theorem}{Theorem}
\newproclaim{definition}{Definition}
\newtheorem{proposition}{Proposition}
\newtheorem{lemma}{Lemma}

\newcommand{\smallset}{A}
\newcommand{\alternateset}{C}
\newcommand{\largeset}{B}
\newcommand{\newset}{\largeset\setminus\smallset}
\newcommand{\setcollection}{\mathcal{A}}
\newcommand{\mle}{\widehat{\theta}}
\newproclaim{remark}{Remark}
\makeatother

\begin{document}
\begin{frontmatter}

\title{Consistency under sampling of exponential random graph models}
\runtitle{Consistency under sampling of exponential graphs}

\begin{aug}
\author[A]{\fnms{Cosma Rohilla} \snm{Shalizi}\ead[label=e1]{cshalizi@cmu.edu}\thanksref{t1}}
\and
\author[A]{\fnms{Alessandro} \snm{Rinaldo}\corref{}\ead[label=e2]{arinaldo@cmu.edu}}
\runauthor{C. R. Shalizi and A. Rinaldo}
\thankstext{t1}{Supported by grants from
the National Institutes of Health (\# 2 R01 NS047493) and the Institute
for New Economic Thinking.}
\affiliation{Carnegie Mellon University}
\address[A]{Department of Statistics\\
Carnegie Mellon University\\
Pittsburgh, PA 15213\\
USA\\
\printead{e1}\\
\phantom{E-mail:\ }\printead*{e2}} %adresu isvedimo komanda gale!
\end{aug}

% HISTORY:
\received{\smonth{11} \syear{2011}}
\revised{\smonth{8} \syear{2012}}

% ABSTRACT

\begin{abstract}
The growing availability of network data and of scientific interest in
distributed systems has led to the rapid development of statistical
models of
network structure. Typically, however, these are models for the entire
network, while the data consists only of a sampled sub-network. Parameters
for the whole network, which is what is of interest, are estimated by
applying the model to the sub-network. This assumes that the model is
\emph{consistent under sampling}, or, in terms of the theory of stochastic
processes, that it defines a projective family. Focusing on the
popular class of exponential random graph models (ERGMs), we
show that this apparently trivial condition is in fact violated by many
popular and scientifically appealing models, and that satisfying it
drastically limits ERGM's expressive power. These results are actually
special cases of more general results about exponential families of dependent
random variables, which we also prove. Using such results, we offer
easily checked conditions for the consistency of maximum likelihood
estimation in ERGMs, and discuss some possible constructive responses.
\end{abstract}

% KEYWORDS
\begin{keyword}[class=AMS]
\kwd[Primary ]{91D30}
\kwd{62B05}
\kwd[; secondary ]{60G51}
\kwd{62M99}
\kwd{62M09}
\end{keyword}
\begin{keyword}
\kwd{Exponential family}
\kwd{projective family}
\kwd{network models}
\kwd{exponential random graph model}
\kwd{sufficient statistics}
\kwd{independent increments}
\kwd{network sampling}
\end{keyword}

\end{frontmatter}

%s1 #&#
\section{Introduction}\label{sec1}

In recent years, the rapid increase in both the availability of data on
networks (of all kinds, but especially social ones) and the demand,
from many
scientific areas, for analyzing such data has resulted in a surge of generative
and descriptive models for network data
\cite{Easley-Kleinberg-networks-crowds-and-markets,MEJN-on-networks}. Within
statistics, this trend has led to a renewed interest in developing, analyzing
and validating statistical models for networks
\cite{Goldenberg-et-al-on-networks,Kolaczyk-on-network-data}. Yet as networks
are a nonstandard type of data, many basic properties of statistical models
for networks are still unknown or have not been properly explored.

In this article we investigate the conditions under which statistical
inferences drawn over a sub-network will generalize to the entire
network. It
is quite rare for the data to ever actually be the \emph{whole} network of
relations among a given set of nodes or units;\setcounter{footnote}{1}\footnote{This sense of the
``whole network'' should not be confused with the technical term ``complete
graph,'' where every vertex has a direct edge to every other vertex.}
typically, only a sub-network is available. Guided by experience of more
conventional problems like regression, analysts have generally fit
models to
the available sub-network, and then extrapolated them to the larger true
network which is of actual scientific interest, presuming that the
models are,
as it were, consistent under sampling. What we show is that this is
only valid
for very special model specifications, and the specifications where it
is \emph{not} valid include some of which are currently among the most popular and
scientifically appealing.

In particular, we restrict ourselves to exponential random graph models
(ERGMs), undoubtedly one of the most important and popular classes of
statistical models of network structure. In addition to the general works
already cited, the reader is referred to
\cite{Frank-Strauss86,Wasserman-Pattison-1996,pstar-primer,Snijders-Pattison-Robins-Handcock-new-specification,Robins-et-al-recent-developments-in-pstar,Wasserman-Robins-dependence-graphs-and-pstar,statnet-special-issue,Park-MEJN-stat-mech-of-networks}
for detailed accounts of these models. There are many reasons ERGMs are so
prominent. On the one hand, ERGMs, as the name suggests, are exponential
families, and so they inherit all the familiar virtues of exponential families
in general: they are analytically and inferentially convenient
\cite{Brown-on-exponential-families}; they naturally arise from considerations
of maximum entropy~\cite{Mandelbrot-sufficiency-and-estimation-in-thermo} and
minimum description length~\cite{Grunwald-on-MDL}, and from
physically-motivated large deviations principles
\cite{Touchette-large-dev-and-stat-mech}; and if a generative model obeys
reasonable-seeming regularity conditions while still having a
finite-dimensional sufficient statistic, it must be an exponential family
\cite{Lauritzen-extremal-families-and-suff-stats}.\footnote{\cite{Mandelbrot-sufficiency-and-estimation-in-thermo}
is still one of the best discussions of the interplay between the formal,
statistical and substantive motivations for using exponential
families.} On
the other hand, ERGMs have particular virtues as models of networks. The
sufficient statistics in these models typically count the number or
density of
certain ``motifs'' or small sub-graphs, such as edges themselves, triangles,
$k$-cliques, stars, etc. These in turn are plausibly related to different
network-growth mechanisms, giving them a substantive interpretation; see,
for example,~\cite{Goodereau-Kitts-Morris-birds-of-a-feather} as an exemplary
application of this idea, or, more briefly, Section~\ref{secergms} below.
Moreover, the important task of edge prediction is easily handled in this
framework, reducing to a conditional logistic regression
\cite{statnet-special-issue}. Since the development of (comparatively)
computationally-efficient maximum-likelihood estimators (based on Monte Carlo
sampling), ERGMs have emerged as flexible and persuasive tools for modeling
network data~\cite{statnet-special-issue}.

Despite all these strengths, however, ERGMs are tools with a serious weakness.
As we mentioned, it is very rare to ever observe the whole network of interest.
The usual procedure, then, is to fit ERGMs (by maximum likelihood or
pseudo-likelihood) to the observed sub-network, and then extrapolate
the same
model, with the same parameters, to the whole network; often this takes the
form of interpreting the parameters as ``provid[ing] information about the
presence of structural effects observed in the network''~\cite{Robins-et-al-recent-developments-in-pstar}, page~194, or the strength of different
network-formation mechanisms;
\cite{Ackland-ONeil-online-collective-identity,Daraganova-et-al-networks-and-geography,de-la-Haye-et-al-on-obesity-and-friendship-networks,Gondal-ERGM-for-knowledge-production,Gonzalez-Bailon-ergm-for-WWW,Schaefer-youth-co-offending-networks,Vermeij-et-al-ergms-of-segregation-networks} are just a few of the more
recent papers doing this. This obviously raises the question of the
statistical (i.e., large sample) consistency of maximum likelihood estimation
in this context. Unnoticed, however, is the logically prior question of
whether it is \emph{probabilistically} consistent to apply the same
ERGM, with
the same parameters, both to the whole network and its sub-networks.
That is,
whether the marginal distribution of a sub-network will be consistent
with the
distribution of the whole network, for all possible values of the model
parameters. The same question arises when parameters are compared across
networks of different sizes (as in, e.g.,
\cite{Faust-Skvoretz-comparing-networks,Goodereau-Kitts-Morris-birds-of-a-feather,Lubbers-Snijders-comparison-of-ergms}).
When this form of consistency fails, then the parameter estimates
obtained from
a sub-network may not provide reliable estimates of, or may not even be
relatable to, the parameters of the whole network, rendering the task of
statistical inference based on a sub-network ill-posed. We formalize this
question using the notion of ``projective families'' from the theory of
stochastic processes. We say that a model is \emph{projective} when the same
parameters can be used for both the whole network and any of its sub-networks.
In this article, we fully characterize projectibility of discrete exponential
families and, as corollary, show that ERGMs are projective only for very
special choices of the sufficient statistic.

%pa1.subsection.subsubsection.1 #&#
\textit{Outline}. Our results are not specific just to
networks, but pertain more generally with exponential families of stochastic
processes. In Section~\ref{secprojectiblity-and-exp-fam}, therefore,
we lay out the
necessary background about projective families of distributions, projective
parameters and exponential families in a somewhat more abstract setting than
that of networks. In Section \ref
{secprojectibility-and-independent-increments} we show
that a necessary and sufficient condition for an exponential family to be
projective is that the sufficient statistics obey a kind of additive
decomposition. This in turn implies strong independence properties. We also
prove results about the consistency of maximum likelihood parameter estimation
under these conditions (Section~\ref{secconsistency-of-mle}). In Section
\ref{secergms}, we apply these results to ERGMs, showing that most popular
specifications for social networks and other stochastic graphs cannot be
projective. We then conclude with some discussion on possible constructive
responses. The proofs are contained in the \hyperref[app]{Appendix}.

%pa1.subsection.subsubsection.2 #&#
\textit{Related work}. An early recognition of the fact that
sub-networks may
have statistical properties which differ radically from those of the whole
network came in the context of studying networks with power-law
(``scale-free'') degree distributions. On the one hand, Stumpf, Wiuf
and May
\cite{Stumpf-Wiuf-May-subnets-are-not-scale-free} showed that
``subnets of
scale-free networks are not scale-free;'' on the other, Achlioptas et al.
\cite{Achlioptas-et-al-bias-of-traceroute} demonstrated that a particular,
highly popular sampling scheme creates the appearance of a power-law degree
distribution on nearly any network. While the importance of network sampling
schemes has been recognized since then~\cite{Kolaczyk-on-network-data}, Chapter 5,
and valuable contributions have come from, for example,
\cite{Kossinets-effects-of-missing-data,Handcock-Gile-sampled-networks,Krivitsky-Handcock-Morris-adjusting-for-network-size,Ahmed-Neville-Kompella-network-sampling},
we are not aware of any work which has addressed the specific issue of
consistency under projection which we tackle here. Perhaps the closest
approaches to our perspective are~\cite{Orbanz-projective-limit-techniques} and
\cite{Xiang-Neville-learning-with-one-network}. The former considers
conditions under which infinite-dimensional families of distributions on
abstract spaces have projective limits. The latter, more concretely, addresses
the consistency of maximum likelihood estimators for exponential
families of
dependent variables, but under assumptions (regarding Markov
properties, the
``shape'' of neighborhoods, and decay of correlations in potential functions)
which are basically incomparable in strength to ours.\looseness=-1

%s2 #&#
\section{Projective statistical models and exponential families}
\label{secprojectiblity-and-exp-fam}

Our results about exponential random graph models are actually special
cases of
more general results about exponential families of dependent random variables,
and are just as easy to state and prove in the general context as for graphs.
Setting this up, however, requires some preliminary definitions and notation,
which make precise the idea of ``seeing more data from the same
source.'' In order to dispense with any measurability issues,
we will implicitly assume the existence of an underlying probability
measure for which the random variables under study are all measurable.
%%Since we will concern ourselves with discrete exponential families of
%distributions, we can safely dispense ourselves from any measurability
%issues.
Furthermore, for the sake of readability we will not rely on the
measure theoretic notion of filtration: though technically appropriate,
it will add nothing to our results.

Let $\setcollection$ be a collection of finite subsets of a
denumerable set
$\indexset$ partially ordered with respect to subset inclusion. For
technical reasons, we will further assume that $\mathcal{A}$ has the
property of being an ideal: that is, if $A$ belongs to $\mathcal{A}$,
then all subsets of $A$ are also in $\mathcal{A}$ and if $A$, and $B$
belongs to $\mathcal{A}$, then so does their union.
We may think of passing from $\smallset$ to $\largeset\supset
\smallset$ as taking increasingly large
samples from a population, or recording increasingly long time series, or
mapping data from increasing large spatial regions, or over an increasingly
dense spatial grid, or looking at larger and larger sub-graphs from a single
network.
Accordingly, we consider the associated
collection of parametric statistical models $\{ \mathcal{P}_{\smallset,\Theta}
\}_{\smallset\in\setcollection}$ indexed by $\setcollection$,
where, for each
$\smallset\in\setcollection$, $\mathcal{P}_{\smallset,\Theta}
\equiv\{
\mathbb{P}_{\smallset,\theta} \}_{\theta\in\Theta}$ is a family of
probability distributions indexed by points $\theta$ in a fixed open set
$\Theta\subseteq\mathbb{R}^d$. The probability distributions in
$\mathcal{P}_{\smallset,\Theta}$ are also assumed to be supported
over the same $\mathcal{X}_{\smallset}$, which are countable\footnote
{Our results
extend to continuous observations straightforwardly, but with annoying
notational overhead.} sets for each $A$. %subset of $\mathbb{R}^{d_{%\smallset}}$.
We assume that the partial order of $\setcollection$ is
isomorphic to the partial order over $\{ \mathcal{X}_{\smallset}
\}_{\smallset\in\setcollection}$, in the sense that $\smallset
\subset
{\largeset}$ if and only if $\mathcal{X}_{\largeset} = \mathcal
{X}_{\smallset}
\times\mathcal{X}_{\newset} $.\vadjust{\goodbreak}

For given $\theta$ and $\smallset$, we
denote with $X_{\smallset}$ the random variable distributed as
$\mathbb{P}_{\smallset,\theta}$. In particular, for a given $\theta
\in\Theta$, we can regard the
$\{ \mathbb{P}_{\smallset,\theta} \}_{\smallset\in\setcollection}$
as finite-dimensional (i.e., marginal) distributions. % of the
%stochastic process $\{ X_A\}$ indexed by $\mathcal{A}$.

\begin{figure}

\includegraphics{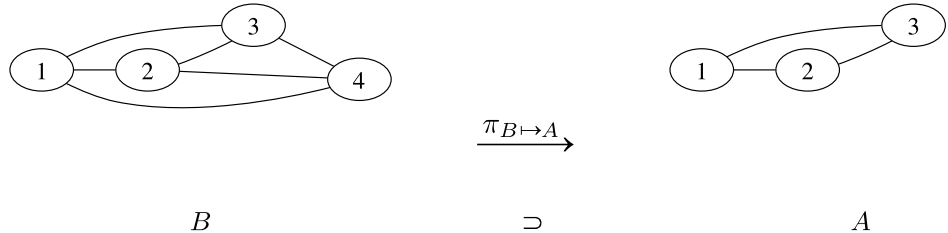}

\caption{Projective structure for networks: when the set of observables
$\smallset$ is contained in the larger set of observables $\largeset$,
$X_{\smallset}$ (on the right) can be recovered from $X_{\largeset}$ (on
the left) through the projection $\pi_{\largeset\mapsto\smallset}$, which
simply drops the extra data.}
\label{figprojective-structure-illustrated}
\end{figure}

For each pair $\smallset, \largeset$ in $\setcollection$ with
$\smallset
\subset\largeset$, we let $\pi_{\largeset\mapsto\smallset} \colon
\mathcal{X}_{\largeset} \rightarrow\mathcal{X}_\smallset$ be the
natural index
projection given by $\pi_{\largeset\mapsto\smallset}(x_{\smallset
},x_{\newset}) = x_{\smallset}$.
In the context of networks, we may think of $\mathbb{I}$ as the set of
nodes of a possibly infinite random graph, which without loss of
generality can be taken to be $\{1,2,\ldots\}$ and of $\mathcal{A}$
as the collection of all finite subsets of $\mathbb{I}$. Then, for
some positive integers $n$ and $m$, we may, for instance, take $A = \{
1,\ldots,n\}$ and $B = \{ 1,\ldots,n,\ldots,n+m\}$, so that
$X_{\smallset}$ will be the induced sub-graph on the first $n$ nodes
and $X_B$ the induced sub-graph on the first $n+m$ nodes.
%$X_{\smallset}$ should be taken as the
%induced sub-graph among the first $n$ (say) nodes within some total
%random graph, and
%$X_{\largeset}$ the sub-graph among the first $n+m$ nodes, where $m >
%0$.
The projection $\pi_{\largeset\mapsto\smallset}$ then just picks
out the
appropriate sub-graph from the larger graph; see Figure
\ref{figprojective-structure-illustrated} for a schematic example. We
will be concerned with a
natural form of probabilistic consistency of the collection $\{
\mathcal{P}_{\smallset,\Theta} \}_{\smallset\in\setcollection}$
which we call
\textit{projectibility}, defined below.

%f1 #&#

% \begin{center}
% \begin{tabular}{ccc}
% \includegraphics[width=0.45\textwidth]{big-graph} & {\Large$
% $\largeset$ & $\supset$ & $\smallset$
% \end{tabular}
% \end{center}
%

%de1 #&#
\begin{definition}
The family $\{ \mathcal{P}_{\smallset,\Theta} \}_{\smallset\in
\setcollection}$ is \emph{projective} if, for any $\smallset$ and $B$
in $\mathcal{A}$ with $\smallset\subset
\largeset$,
%
%e1 #&#
\begin{equation}
\mathbb{P}_{\smallset,\theta} = \mathbb{P}_{\largeset,\theta} \circ\pi_{\largeset\mapsto\smallset}^{-1}\qquad
\forall\theta \in\Theta.
\end{equation}
\end{definition}

See~\cite{Kallenberg-mod-prob}, page 115, for more general treatment of
projectibility. In words, $\{ \mathcal{P}_{\smallset,\Theta} \}
_{\smallset\in
\setcollection}$ is a projective family when $\smallset\subset
\largeset$
implies that $\mathbb{P}_{\smallset,\theta} $ can be recovered by
marginalization over $\mathbb{P}_{\largeset,\theta}$, for all
$\theta$. % (Figure~\ref{figprojective-family}).
Within a projective family,
$\mathbb{P}_{\theta}$ denotes the infinite-dimensional distribution,
which thus
exists by the Kolmogorov extension theorem~\cite{Kallenberg-mod-prob}, Theorem 6.16, page 115.

Projectibility is automatic when the generative model calls for
independent and
identically distributed (IID) observations. It is also generally unproblematic
when the model is specified in terms of \emph{conditional}
distributions: one
then just uses the Ionescu Tulcea extension theorem in place of that of
Kolmogorov~\cite{Kallenberg-mod-prob}, Theorem~6.17, page~116.
However, many
models are specified in terms of \emph{joint} distributions for various index
sets, and this, as we show in Theorem \ref
{thmprojectible-iff-separable}, can
rule out projectibility.

We restrict ourselves to \emph{exponential family} models by assuming
that, for
each choice of $\theta\in\Theta$ and $\smallset\in\setcollection$,
$\mathbb{P}_{\smallset,\theta}$ has density with respect to the counting
measure over $\mathcal{X}_{\smallset}$ given by
%
%e2 #&#
\begin{equation}
p_{\smallset,\theta}(x) = \frac{e^{\langle\theta,
t_\smallset(x)\rangle}}{z_\smallset(\theta)},\qquad x \in\mathcal{X}_{\smallset},
\label{eqnexponential-family-density}
\end{equation}
where $t_{\smallset} \colon\mathcal{X}_{\smallset} \rightarrow
\mathbb{R}^d$
is the measurable function of minimal sufficient statistics, and
$z_{\smallset} \colon
\Theta\rightarrow\mathbb{R}$ is the \emph{partition function} given by
%
%e3 #&#
\begin{equation}
z_{\smallset}(\theta) \equiv\sum_{x \in\mathcal{X}_{\smallset}}
e^{\langle\theta, t_\smallset(x)\rangle}.
\end{equation}
If $X_{\smallset} \sim\mathbb{P}_{\smallset,\theta}$, we will write
$T_{\smallset} \equiv t_{\smallset}(X_{\smallset})$ for the random variable
corresponding to the sufficient statistic. Equation  \eqref{eqnexponential-family-density} implies that $T_{\smallset}$
itself has
an exponential family distribution, with the same parameter $\theta$ and
partition function $z_{\smallset}(\theta)$~\cite{Brown-on-exponential-families}, Proposition
1.5. Specifically, the distribution function
is
%
%e4 #&#
\begin{equation}
\mathbb{P}_{\smallset,\theta}(T_{\smallset} = t) = \frac
{e^{\langle\theta, t\rangle} v_{\smallset}(t)}{z_{\smallset
}(\theta)},
\end{equation}
where the term $v_{\smallset}(t) \equiv\llvert  \{ x \in\mathcal
{X}_A \colon t_A(x) = t \}\rrvert $, which we will call the \textit{volume factor},
counts the number of points in $\mathcal{X}_{\smallset}$ with the
same sufficient statistics $t$. The moment generating function of
$T_{\smallset}$ is
%
%e5 #&#
\begin{equation}
M_{\theta,\smallset}(\phi) = \mathbf{E}_{\theta}\bigl[ e^{\langle\phi, T_{\smallset}\rangle} \bigr]
= z_{\smallset}(\theta+\phi)/z_{\smallset}(\theta). \label{eqnmgf-of-exponential-family}
\end{equation}

If the sufficient statistic is completely additive, that is, if
$t_\smallset(x_\smallset) =\break \sum_{i\in\smallset}{t_{\{i\}}(x_i)}$,
then this
is a model of independent (if not necessarily IID) data. In general, however,
the choice of sufficient statistics may impose, or capture, dependence between
observations.

Because we are considering exponential families defined on increasingly large
sets of observations, it is convenient to introduce some notation
related to
multiple statistics. Fix $\smallset, \largeset\in\setcollection$
such that
$\smallset\subset\largeset$. Then $t_{\largeset}\dvtx \mathcal
{X}_{\largeset}
\mapsto\mathbb{R}^d$, and we will sometimes write this function $t(x,y)$,
where the first argument is in $\mathcal{X}_{\smallset}$ and the
second in
$\mathcal{X}_{\newset}$. We will have frequent recourse to the
increment to
the sufficient statistic, $t_{\newset}(x,y) \equiv t_{\largeset}(x,y) -
t_{\smallset}(x)$. The volume factor
$v_{\largeset}(t_{\largeset}(x_{\largeset}))$ is defined as before,
but we shall
also consider, for each observable value $t$ of the sufficient
statistics for
$A$ and increment $\delta$ of the sufficient statistics from $A$ to
$B$, the
\textit{joint volume factor,}
%
%e6 #&#
\begin{equation}
\label{eqjointvolfactor} v_{\smallset,\newset}(t,\delta) \equiv\bigl\llvert \bigl\{
(x,y) \in\mathcal {X}_{\largeset} \colon t_{\smallset}(x) = t \mbox{ and }
t_{\newset
}(x,y) = \delta\bigr\} \bigr\rrvert,
\end{equation}
and the \textit{conditional volume factor,}
%
%e7 #&#
\begin{equation}
v_{\newset|\smallset}(\delta,x) \equiv\bigl\llvert \bigl\{ y \in\mathcal
{X}_{\newset} \colon t_{\newset}(x,y) = \delta\bigr\} \bigr\rrvert.
\end{equation}
As we will see, these volume factors play a key role in characterizing
projectibility.

%s3 #&#
\section{Projective structure in exponential families}
\label{secprojectibility-and-independent-increments}

In this section we characterize projectibility in terms of the
increments of
the vector of sufficient statistics. In particular we show that exponential
families are projective if, and only if, their sufficient statistics decompose
into separate additive contributions from disjoint observations in a
particularly nice way which we formalize in the following
definition.\vspace*{-2pt}

%de2 #&#
\begin{definition} The sufficient statistics of the family $\{
\mathcal{P}_{\smallset,\Theta} \}_{\smallset\in\setcollection}$
have \emph{separable increments} when, for each $\smallset\subset\largeset$, $x
\in
\mathcal{X}_A$, the range of possible increments $\delta$ is the same
for all
$x$, and the conditional volume factor is constant in $x$, that is,
$v_{\newset|\smallset}(\delta,x) = v_{\newset}(\delta)$.\vspace*{-2pt}
\end{definition}

It is worth noting that the property of having separable increments is an
intrinsic property of the family $\{ \mathcal{P}_{\smallset,\Theta}
\}_{\smallset\in\setcollection}$ that depends only on the functional
forms of
the sufficient statistics $\{ t_A\}_{A \in\mathcal{A}}$ and not on
the model
parameters $\theta\in\Theta$. This follows from the fact that, for
any $A$,
the probability distributions $\{ \mathbb{P}_{A,\theta}\}_{\theta\in
\Theta}$
have identical support $\mathcal{X}_A$. Thus, this property holds for
all of
$\theta$ or none of them.

The main result of this paper is then as follows.\vspace*{-2pt}

%th1 #&#
\begin{theorem}
The exponential family $\{ \mathcal{P}_{\smallset,\Theta} \}
_{\smallset\in
\setcollection}$ is projective if and only if the sufficient
statistics $\{
T_A \}_{A \in\mathcal{A}}$ have separable increments.
\label{thmprojectible-iff-separable}\vspace*{-2pt}
\end{theorem}

%s3.1 #&#
\subsection{Independence properties}
\label{secindependence-properties}

Because projectibility implies separable increments, it also carries
statistical-independence implications. Specifically, it implies that the
increments to the sufficient statistics are statistically independent,
and that
$X_{\largeset\setminus\smallset}$ and $X_{\smallset}$ are conditionally
independent given increments to the sufficient statistic. Interestingly,
independent increments for the statistic are necessary but not quite sufficient
for projectibility. These claims are all made more specific in the
propositions which follow.

We first show that projectibility implies that the sufficient
statistics have
independent increments. In fact, a stronger results holds, namely that
the increments of the sufficient
statistics are independent of the actual sequence. Below we will write
$T_{\newset}$ to signify $T_{\largeset} - T_{\smallset}$.\vspace*{-2pt}
%
%pr1 #&#
\begin{proposition}
If the exponential family $\{ \mathcal{P}_{\smallset,\Theta} \}
_{\smallset
\in\setcollection}$ is projective, then sufficient statistics $\{ T_A
\}_{A \in\mathcal{A}}$ have independent increments, that is,
$\smallset\subset
\largeset$ implies that $T_{\largeset} - T_{\smallset} \indep
T_\smallset$
under all $\theta$.
\label{propproj-to-II}\vspace*{-2pt}
\end{proposition}

%pr2 #&#
\begin{proposition}
In a projective exponential family, $T_{\newset}\indep X_{\smallset}$.
\label{propincrement-indep-of-old-observable}\vspace*{-2pt}
\end{proposition}

We note that independent increments for the sufficient statistics
$T_{\smallset}$ in no way implies independence of the actual observations
$X_{\smallset}$. As a simple illustration, take the one-dimensional Ising
model,\footnote{Technically, with ``free'' boundary conditions; see
\cite{Landau-Lifshitz}.} where\vadjust{\goodbreak} $\indexset= \mathbb{N}$, each
${\mathcal{X}}_i = \pm1$, $\setcollection$ consists of all intervals
from $1$
to $n$, and the single sufficient statistic $T_{1:n} = \sum_{i=1}^{n-1}{X_i
X_{i+1}}$. Clearly, $T_{1:(n+1)} - T_{1:n} = +1$ when $X_n = X_{n+1}$,
otherwise $T_{1:(n+1)} - T_{1:n} = -1$. Since $v_{1:(n+1)|1:n}(+1,x) =
v_{1:(n+1)|1:n}(-1,x) = 1$, by Theorem \ref
{thmprojectible-iff-separable}, the
model is projective. By Proposition~\ref{propproj-to-II}, then,
increments of
$T$ should be independent, and direct calculation shows the probability of
increasing the sufficient statistic by 1 is $e^{\theta}/(1+e^{\theta
})$, no
matter what $X_1, \ldots, X_n$ are. While the sufficient statistic has
independent increments, the random variables $X_i$ are all dependent on one
another.\footnote{Note that while this is a \emph{graphical} model, it
is not a
model of a random graph. (The graph is rather the one-dimensional lattice.)
Rather, it is used here merely to exemplify the general result about
exponential families. We turn to exponential random graph models in Section
\ref{secergms}.}
\label{example:ising-model}

The previous results provide a way, and often a simple one, for checking
whether projectibility fails: if the sufficient statistics do not have
independent increments, then the family is not projective. As we will
see, this
test covers many statistical models for networks.

It is natural to inquire into the converse to these propositions. It is fairly
straightforward (if somewhat lengthy) to show that independent
increments for
the sufficient statistics implies that the joint volume factor separates.

%pr3 #&#
\begin{proposition}
If an exponential family has independent increments, $T_{\newset}
\indep
T_{\smallset}$, then its joint volume factor separates,
$v_{\smallset,\newset}(t,\delta) = v_{\smallset}(t)v_{\newset
}(\delta)$, and
the distribution of $T$ is projective.
\label{propII-to-volume-factor-separation}
\end{proposition}

However, independent increments for the sufficient statistics do \emph{not}
imply that separable increments (hence projectibility), as shown by the next
counter-example. Hence independent increments are a necessary but not
sufficient condition for projectibility.

Suppose that $\mathcal{X}_{\smallset} =  \{ a,b,c,d \}$, and
$\mathcal{X}_{\newset} =  \{i, \mathit{ii}, \mathit{iii}, {iv}, v \}$. (Thus
there are 20
possible values for $X_{\largeset}$.) Let
\begin{eqnarray*}
+1 & = & t_{\smallset}(a) = t_{\smallset}(b),
\\
-1 & = & t_{\smallset}(c) = t_{\smallset}(d)
\end{eqnarray*}
so that $v_{\smallset}(+1) = v_{\smallset}(-1) = 2$. Further, let
\begin{eqnarray*}
2 & = & t_{\largeset}(a,i) = t_{\largeset}(a,\mathit{ii}),
\\
0 & = & t_{\largeset}(a,\mathit{iii}) = t_{\largeset}(a,iv) = t_{\largeset
}(a,v),
\\
0 & = & t_{\largeset}(b,i) = t_{\largeset}(b,\mathit{ii}),
\\
2 & = & t_{\largeset}(b,\mathit{iii}) = t_{\largeset}(b,iv) = t_{\largeset
}(b,v),
\\
t_{\largeset}(c,y) & = & t_{\largeset}(a,y) -2,
\\
t_{\largeset}(d,y) & = & t_{\largeset}(b,y) -2.
\end{eqnarray*}
It is not hard to verify that $T_{\newset}$ is always either $+1$ or
$-1$. It
is also straightforward to check that $v_{\smallset,\newset}(t,\delta
) = 5$ for
all combinations of $t$ and $\delta$, implying that $v_{\newset}(+1) =
v_{\newset}(-1) = 2.5$, and that the joint volume factor separates. On the
other hand, the \emph{conditional} volume factors are not constant in
$x$, as
$v_{\newset|\smallset}(+1,a)=2$ while $v_{\newset|\smallset}(+1,b)
= 3$. Thus,
the sufficient statistic has independent increments, but does not have
separable increments. Since projective families have separable increments
(Proposition~\ref{propproj-to-SI}), this cannot be a projective
family. (This
can also be checked by a direct and straightforward, if even more tedious,
calculation.)

We conclude this section with a final observation. Butler, in
\cite{Butler-predictive-likelihood}, showed that when observations
follow from
an IID model with a minimal sufficient statistic, the predictive distribution
for the next observation can be written entirely in terms of how different
hypothetical values would change the sufficient statistic; cf.
\cite{Lauritzen-sufficiency-and-prediction,Besag-candidates-formula}. This
predictive sufficiency property carries over to our setting.

%th2 #&#
\begin{theorem}[(Predictive sufficiency)]\label{thmpredictive-sufficiency}
In a projective exponential family, the distribution of $X_{\newset}$
conditional on $X_{\smallset}$ depends on the data only through~$T_{\newset}$.
\end{theorem}

\begin{figure}

\includegraphics{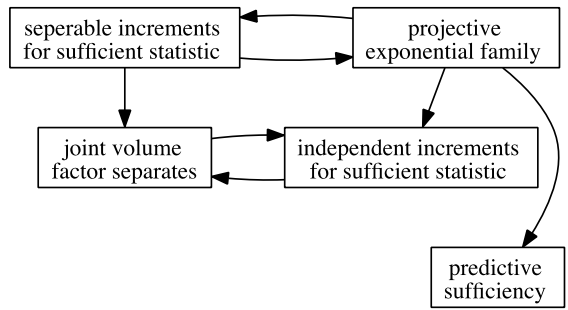}

\caption{Relations among the main properties of models considered in Section
\protect\ref{secprojectibility-and-independent-increments}. Probabilistic
properties of the models are on the right, and algebraic/combinatorial
properties
of the sufficient statistic are on the left.}
\label{figimplications}
\end{figure}

The main implications among our results are summarized in Figure
\ref{figimplications}.

%f2 #&#

%s3.2 #&#
\subsection{Remarks, applications and extensions}
\label{secgrowing-parameter-set}
\mbox{}

\textit{Exponential families of time series}. As the example of the Ising
model in Section~\ref{secindependence-properties} (page
\pageref{example:ising-model}) makes clear, our theorem applies
whenever we
need an exponential family to be projective, not just when the data are
networks. In particular, they apply to exponential families of time series,
where $\indexset$ is the natural or real number line (or perhaps just its
positive part), and the elements of $\setcollection$ are intervals. An
exponential family of stochastic processes on such a space has projective
parameters if, and only if, its sufficient statistics have separable
increments, and so only if they have independent increments.

\textit{Transformation of parameters}. Allowing the dimension of
$\theta$ to
be fixed, but for its components to change along with $\smallset$,
does not
really get out of these results. Specifically, if $\theta$ is to be re-scaled
in a way that is a function of $\smallset$ alone, we can recover the
case of a
fixed $\theta$ by ``moving the scaling across the inner product,''
that is, by
re-defining $T_{\smallset}$ to incorporate the scaling. With a
sample-invariant $\theta$, it is this transformed $T$ which must have separable
increments. Other transformations can either be dealt with similarly, or
amount to using a nonuniform base measure; see below.

\textit{Statistical-mechanical interpretation}. It is interesting to consider
the interpretation of our theorem, and of its proof, in terms of statistical
mechanics. As is well known, the ``canonical'' distributions in statistical
mechanics are exponential families (Boltzmann--Gibbs distributions),
where the
sufficient statistics are ``extensive'' physical observables, such as energy,
volume, the number of molecules of various species, etc., and the natural
parameters are the corresponding conjugate ``intensive'' variables,
such as,
respectively, (inverse) temperature, pressure, chemical potential,
etc.~\cite{Landau-Lifshitz,Mandelbrot-sufficiency-and-estimation-in-thermo}.
Equilibrium between two systems, which interact by exchanging the variables
tracked by the extensive variables, is obtained if and only if they
have the same
values of the intensive parameters~\cite{Landau-Lifshitz}. In our
terms, of
course, this is simply projectibility, the requirement that the same parameters
hold for all sub-systems. What we have shown is that for this to be
true, the
increments to the extensive variables must be completely unpredictable from
their values on the sub-system.

Furthermore, notice the important role played in both halves of the
proof by
the separation of the joint volume factor, $v_{\smallset,\newset
}(t,\delta) =
v_{\smallset}(t) v_{\newset}(\delta)$. In terms of statistical
mechanics, a
macroscopic state is a collection of microscopic configurations with
the same
value of one or more macroscopic observables. The Boltzmann entropy of a
macroscopic state is (proportional to) the logarithm of the volume of those
microscopic states~\cite{Landau-Lifshitz}. If we define our
macroscopic states
through the sufficient statistics, then their Boltzmann entropy is just
$\log{v}$. Thus, the separation of the volume factor is the same as the
additivity of the entropy across different parts of the system, that
is, the
entropy is ``extensive.'' Our results may thus be relevant to debates in
statistical mechanics about the appropriateness of alternative, nonextensive
entropies; cf.~\cite{Nauenberg-critique-of-q-entropy}.

\textit{Beyond exponential families}. It is not clear just how
important it
is that we have an exponential family, as opposed to a family admitting a
finite-dimensional sufficient statistic. As is well known, the two concepts
coincide under some regularity conditions
\cite{Barndorff-Nielsen-exponential-families}, but not quite strictly,
and it
would be interesting to know whether or not the exponential form of equation
\eqref{eqnexponential-family-density} is strictly required. We have attempted
to write the proofs in a way which minimizes the use of this form (in
favor of
the Neyman factorization, which only uses sufficiency), but have not succeeded
in eliminating it completely. We return to this matter in the conclusion.

\textit{Prediction}. We have focused on the implications of projectibility
for parametric inference. Exponential families are, however, often used in
statistics and machine learning as generative models in applications
where the
only goal is prediction~\cite{Wainwright-Jordan-graphical-models-and-exponential-families},
and so (to
quote Butler~\cite{Butler-predictive-likelihood}) ``all parameters are nuisance
parameters.'' But even in then, it must be possible to consistently
extend the
generative model's distribution for the training data to a joint distribution
for training and testing data, with a single set of parameters shared
by both
old and new data. While this requirement may seem too trivial to
mention, it
is, precisely, projectibility.

\textit{Growing number of parameters}. In the proof of Theorem
\ref{thmprojectible-iff-separable}, we used the fact that
$T_{\smallset}$, and
hence $\theta$, has the same dimension for all $\smallset\in
\setcollection$.
There are, however, important classes of models where the number of parameters
is allowed to grow with the size of the sample. Particularly important, for
networks, are models where each node is allowed a parameter (or two) of its
own, such as its expected degree; see, for instance, the classic $p_1$ model
of~\cite{Holland-Leinhardt-p1}, or the ``degree-corrected block
models'' of
\cite{Karrer-MEJN-blockmodels-and-community-structure}. We can
formally extend
Theorem~\ref{thmprojectible-iff-separable} to cover some of these
cases---including those two particular specifications---as follows.

Assume that $T_{\smallset}$ has a dimension which is strictly
nondecreasing as
$\smallset$ grows, that is, $d_{\smallset} \leq d_{\largeset}$
whenever $\smallset
\subset\largeset$. Furthermore, assume that the set of parameters
$\theta_{\smallset}$ only grows, and that the meaning of the old
parameters is
not disturbed. That is, under projectibility we should have
%
%e8 #&#
\begin{equation}
\mathbb{P}_{\largeset,\theta_{\largeset}} \cdot \pi_{\largeset\mapsto\smallset}^{-1} =
\mathbb{P}_{\smallset,\pi_{\mathbb{R}^{d_\largeset}\mapsto\mathbb
{R}^{d_\smallset}}
\theta_{\largeset}}(\cdot).
\end{equation}
For any fixed pair $\smallset\subset\largeset$, we can accommodate this
within the proof of Theorem~\ref{thmprojectible-iff-separable} by re-defining
$T_{\smallset}$ to be a mapping from $\mathcal{X}_{\smallset}$ to
$\mathbb{R}^{d_{\largeset}}$, where the extra $d_{\largeset} -
d_{\smallset}$
components of the vector are always zero. The extra parameters in
$\theta_{\largeset}$ then have no influence on the distribution of
$X_{\smallset}$ and are unidentified on $\smallset$, but we have, formally,
restored the fixed-parameter case. The ``increments'' of the extra components
of $T_{\largeset}$ are then simply their values on $X_{\largeset}$,
and, by the
theorem, the range of values for these statistics, and the number of
configurations on $\mathcal{X}_{\newset}$ leading to each value, must
be equal
for all $x \in\mathcal{X}_{\smallset}$.

Adapting our conditions for the asymptotic convergence of maximum likelihood
estimators (Section~\ref{secconsistency-of-mle}) to the
growing-parameter setting
is beyond our scope here.

\textit{Nonuniform base measures}.
If the exponential densities in \eqref{eqnexponential-family-density} are
defined with respect to nonuniform base measures different from the counting
measures, the sufficient statistics need not have separable increments.
In the supplementary material~\cite{supp}
we address this issue and describe the modifications and additional
assumptions required for our analysis to remain valid. We thank an anonymous
referee and Pavel Krivitsky for independently brining up this subtle
point to
our attention.

%s4 #&#
\section{Consistency of maximum likelihood estimators}
\label{secconsistency-of-mle}

Statistical inference in an exponential family naturally centers on the
parameter $\theta$. As is well known, the maximum likelihood estimator
$\mle$
takes a particularly simple form, obtainable using the fact [which
follows from
equation \eqref{eqnmgf-of-exponential-family}] that $\nabla_{\theta}
z_{\smallset}(\theta) =
z_{\smallset}(\theta)\mathbf{E}_{\theta}[ T_{\smallset} ]$,
%
%e9 #&#
\begin{eqnarray}\label{eqnimplicit-exp-family-mle}
0 =  \nabla_{\theta} \frac{e^{\langle\theta, t_{\smallset
}(x)\rangle}}{z_{\smallset}(\theta)}\Big| _{\theta=\mle} & = &
\frac{-z_{\smallset}(\mle) t_{\smallset}(x)e^{\langle\mle,
t_{\smallset}(x_{\smallset})\rangle} + e^{\langle\mle,
t_{\smallset}(x)\rangle} z_{\smallset}(\mle)\mathbf{E}_{\mle}[
T_{\smallset} ]}{z_{\smallset}^2(\theta)},
\nonumber
\\[-8pt]
\\[-8pt]
\nonumber
t_{\smallset}(x) & = & \mathbf{E}_{\mle}[ T_{\smallset} ].
\end{eqnarray}
In words, the most likely value of the parameter is the one where the expected
value of the sufficient statistic equals the observed value.

Assume the conditions of Theorem~\ref{thmprojectible-iff-separable}
hold, so
that the parameters are projective and the sufficient statistics have (by
Lemma~\ref{lemmaseparable-volume-factor-implies-projectible-increments})
independent increments. Define the logarithm of the partition function
$a_{\smallset}(\theta) \equiv\log{z_{\smallset}(\theta
)}$.\footnote{In
statistical mechanics, $-a_{\smallset}$ would be the Helmholtz free energy.}
Suppose that
%
%e10 #&#
\begin{equation}
a_{\smallset}(\theta) = r_{|\smallset|} a(\theta), \label
{eqnscaling-factor}
\end{equation}
where $|\smallset|$ is some positive-valued measure of the size of
$\smallset$,
$r_{|\smallset|}$ a positive monotone-increasing function of it and $a
\colon
\Theta\mapsto\mathbb{R}$ is differentiable (at least at~$\theta$).
Then, by
equation \eqref{eqnmgf-of-exponential-family} for the moment
generating function,
the cumulant generating function of $T_{\smallset}$ is
%
%e11 #&#
\begin{equation}
\kappa_{\smallset,\theta}(\phi) = r_{|\smallset|}\bigl(a(\theta+\phi) - a(\theta)
\bigr).
\end{equation}
From the basic properties of cumulant generating functions, we have
%
%e12 #&#
\begin{equation}
\mathbf{E}_{\theta}[ T_{\smallset} ] = \nabla_{\phi} \kappa
_{\smallset,\theta}(0) = r_{|\smallset|} \nabla a(\theta). \label{eqnexpect-of-suff-stat}
\end{equation}
Substituting into equation \eqref{eqnimplicit-exp-family-mle},
%
%e13 #&#
\begin{equation}
\frac{t_{\smallset}(x)}{r_{|\smallset|}} = \nabla a(\mle).
\end{equation}
Thus to control the convergence of $\mle$, we must control the
convergence of
$T_{\smallset}/r_{|\smallset|}$.

Consider a growing sequence of sets $\smallset$ such that
$r_{|\smallset|}
\rightarrow\infty$. Since $T_{\smallset}$ has independent
increments, and the
cumulant generating functions for different $\smallset$ are all\vadjust{\goodbreak}
proportional to
each other, we may regard $T_{\smallset}$ as a time-transformation of a
L{\'e}vy process $Y_r$~\cite{Kallenberg-mod-prob}. That is, there is a
continuous-time stochastic process $Y$ with IID increments, such that
$Y_1$ has
cumulant generating function $a(\theta+\phi) - a(\theta)$, and
$T_{\smallset} =
Y_{r_{|\smallset|}}$. Note that $T_{\smallset}$ itself does not have
to have
IID increments, but rather the distribution of the increment
$T_{\largeset} -
T_{\smallset}$ must only depend on $r_{|\largeset|} - r_{|\smallset|}$.
Specifically, from Lemma \ref
{corconditional-moment-generating-function} and
equation \eqref{eqnscaling-factor}, the cumulant generating function
of the
increment must be $(r_{|\largeset|} -
r_{|\smallset|})[a(\theta+\phi)-a(\theta)]$. The scaling factor homogenizes
(so to speak) the increments of $T$.

Writing the sufficient statistic as a transformed L{\'e}vy process
yields a
simple proof that $\mle$ is strongly (i.e., almost-surely) consistent.
Since a
L{\'e}vy process has IID increments, by the strong law of large numbers
$Y_{r_{|\smallset|}}/r_{|\smallset|}$ converges almost surely
($\mathbb{P}_{\theta}$) to $\mathbf{E}_{\theta}[ Y_1 ]$
\cite{Kallenberg-mod-prob}. Since $T_{\smallset} = Y_{r_{|\smallset
|}}$, it
follows that $T_{\smallset}/r_{|\smallset|} \rightarrow
\mathbf{E}_{\theta}[ Y_1 ]$ a.s. ($\mathbb{P}_{\theta}$) as well; but
this limit
is $\nabla a(\theta)$. Thus the MLE converges on $\theta$ almost
surely. We
have thus proved

%th3 #&#
\begin{theorem}
Suppose that the model $\mathbb{P}_{\theta}$ is projective, and that
the log
partition function obeys equation \eqref{eqnscaling-factor} for each
$\smallset\in\setcollection$. Then the maximum likelihood estimator exists
and is strongly consistent.
\end{theorem}

We may extend this in a number of ways. First, if the scaling relation
equation~\eqref{eqnscaling-factor} holds for a particular $\theta$ (or set of
$\theta$),
then $T_{\smallset}/r_{|\smallset|}$ will converge almost surely for that
$\theta$. Thus, strong consistency of the MLE may in fact hold over certain
parameter regions but not others. Second, when $d > 1$, all components of
$T_{\smallset}$ must be scaled by the \emph{same} factor $r_{|\smallset|}$.
Making the expectation value of one component of $T$ be $O(|\smallset
|)$ while
another was $O(|\smallset|^3)$ (e.g.) would violate equation
\eqref{eqnexpect-of-suff-stat} and so equation \eqref
{eqnscaling-factor} as
well.

Finally, while the exact scaling of equation \eqref{eqnscaling-factor},
together with the independence of the increments, leads to strong consistency
of the MLE, ordinary consistency (convergence in probability) holds under
weaker conditions. Specifically, suppose that log partition function or free
energy scales in the limit as the size of the assemblage grows,
%
%e14 #&#
\begin{equation}
\lim_{r_{|\smallset|} \rightarrow\infty}{a_{\smallset}(\theta )/r_{|\smallset|}} = a(
\theta); \label{eqnasymptotic-scaling}
\end{equation}
we give examples toward the end of Section~\ref{secergms} below. We may then
use the following theorem:

%th4 #&#
\begin{theorem}
Suppose that an exponential family shows approximate scaling, that is, equation
\eqref{eqnasymptotic-scaling} holds, for some $\theta$. Then, for any
measurable set $K \subseteq\mathbb{R}^d$,
%
%e15 #&#
%e16 #&#
\begin{eqnarray}
\label{eqnldp-lower} \liminf_{r_{|\smallset|} \rightarrow\infty}{\frac{1}{r_{|\smallset
|}}\log{
\mathbb{P}_{\smallset,\theta} \biggl(\frac{T_{\smallset
}}{r_{|\smallset|}} \in K \biggr)}} & \geq& -\inf
_{t \in\mathrm
{int} K}{J(t)},
\\
\limsup_{r_{|\smallset|} \rightarrow\infty}{\frac{1}{r_{|\smallset
|}}\log{\mathbb{P}_{\smallset,\theta}
\biggl(\frac{T_{\smallset
}}{r_{|\smallset|}} \in K \biggr)}} & \leq& -\inf_{t \in\mathrm
{cl} K}{J(t)},
\label{eqnldp-upper}
\end{eqnarray}
where
%
%e17 #&#
\begin{equation}
J(t) = \sup_{\phi\in\mathbb{R}^d}{\langle\phi, t \rangle- \bigl[a(\theta+\phi)
- a(\theta)\bigr]}, \label{eqnrate-function}
\end{equation}
and $\mathrm{int} K$ and $\mathrm{cl} K$ are, respectively, the
interior and
the closure of $K$.
\label{thmldp}
\end{theorem}

When the limits in equations \eqref{eqnldp-lower} and \eqref{eqnldp-upper}
coincide, which they will for most nice sets $K$, we may say that
%
%e18 #&#
\begin{equation}
\frac{1}{r_{|\smallset|}}\log{\mathbb{P}_{\smallset,\theta} \biggl(\frac{T_{\smallset}}{r_{|\smallset|}} \in K
\biggr)} \rightarrow -\inf_{t \in K}{J(t)} \label{eqnldp-for-suff-stat}.
\end{equation}
Since $J(t)$ is minimized at 0 when $t = \nabla a(\theta)$,\footnote
{For small
$\varepsilon\in\mathbb{R}^d$, by a second order Taylor expansion,
$J(\varepsilon
+ \nabla a(\theta)) \approx\frac{1}{2}\langle\varepsilon, I(\theta)
\varepsilon\rangle$, where $I(\theta)$ acts as the Fisher information rate;
cf.~\cite{Bahadur-limit-theorems}.} equation \eqref{eqnldp-for-suff-stat}
holds in particular for any neighborhood of $\nabla a(\theta)$, and
for the
complement of such neighborhoods, where the infimum of $J$ is strictly
positive. Thus $T_{\smallset}/r_{|\smallset|}$ converges in
probability to
$\nabla a(\theta)$, and $\mle\stackrel{P}{\rightarrow} \theta$,
for all
$\theta$ where equaiton \eqref{eqnasymptotic-scaling}
holds.\looseness=-1

Heuristically, when equation \eqref{eqnasymptotic-scaling} holds but equation
\eqref{eqnscaling-factor} fails, we may imagine approximating the actual
collection of dependent and heterogeneous random variables with an
average of
IID, homogenized effective variables, altering the behavior of the global
sufficient statistic $T$ by no more than $o_P(r_{|\smallset|})$. In
statistical-mechanical terms, this means using renormalization~\cite{Yeomans}.
Probabilistically, the existence of a limiting (scaled) cumulant generating
function is a weak dependence condition~\cite{den-Hollander-large-deviations}, Section
V.3.2. While under equation
\eqref{eqnscaling-factor} we identified the $T_{\smallset}$ process
with a
time-transformed L{\'e}vy process, now we can only use a central limit theorem
to say they are close~\cite{den-Hollander-large-deviations}, Section
V.3.1, reducing
almost-sure to stochastic convergence; see~\cite{Jona-Lasinio-RG-and-prob-theory} on the relation between central limit
theorems and renormalization. In any event, asymptotic scaling of the log
partition function implies $\mle$ is consistent.

%s5 #&#
\section{Application: Nonprojectibility of exponential random graph models}
\label{secergms}

As mentioned in the \hyperref[sec1]{Introduction}, our general results about projective
structure in exponential families arose from questions about
exponential random
graph models of networks. To make the application clear, we must fill
in some
details regarding ERGMs.

Given a group of $n$ nodes, the network among them is represented by
the binary
$n \times n$ \emph{adjacency matrix} $X$, where $X_{ij} = 1$ if there is
a tie
from $i$ to $j$ and is 0 otherwise. (Undirected graphs impose $X_{ij} =
X_{ji}$.) We may also have covariates for each node, say $Y_i$. Our
projective structure will in fact be that of looking at the sub-graphs among
larger and larger groups of nodes. That is, $\smallset$ is the sub-network
among the first $n$\vadjust{\goodbreak} nodes, and $\largeset\supset\smallset$ is the sub-network
among the first $n+m$ nodes. The graph or adjacency matrix itself is the
stochastic process which is to have an exponential family distribution,
conditional on the covariates
%
%e19 #&#
\begin{equation}
p_{\theta}(x|y) = \frac{e^{\langle\theta, t(x,y) \rangle
}}{z(\theta|y)}.
\end{equation}
(We are only interested in the exponential-family distribution of the graph
holding the covariates fixed.) As mentioned above, the components of $T$
typically count the number of occurrences of various sub-graphs or
motifs---as edges, triangles, larger cliques, ``$k$-stars'' ($k$ nodes connected through
a central node), etc.---perhaps interacted with values of the nodal
covariates. The definition of $T$ may include normalizing the counts of these
``motifs'' by data-independent combinatorial factors to yield densities.

A \emph{dyad} consists of an unordered pair of individuals. In a dyadic
independence model, each dyad's configuration is independent of every other
dyad's (conditional on $Y$). In an ERGM, dyadic independence is
equivalent to
the (vector-valued) statistic $T$ adding up over dyads,
%
%e20 #&#
\begin{equation}
t(X,Y) = \sum_{i=1}^{n}{\sum
_{j < i}{t_{ij}(X_{ij}, X_{ji},
Y_i, Y_j)}}.
\end{equation}
That is, the statistic can be written as a sum of terms over the information
available for each dyad. In particular, in \emph{block models}
\cite{Bickel-Chen-on-modularity}, $Y_i$ is categorical, giving the
type of node
$i$, and the vector of sufficient statistics counts dyad configurations among
pairs of nodes of given pairs of types. Dyadic independence implies
projectibility: since all dyads have independent configurations, each dyad
makes a separate additive contribution to $T$. Going from $n-1$ to $n$ nodes
thus adds $n$ terms, unconstrained by the configuration among the $n-1$ nodes.
$T$ thus has separable increments, implying projectibility by Theorem
\ref{thmprojectible-iff-separable}. (Adding a new node adds only edges
between the old nodes and the new, without disturbing the old
counts.)\footnote{We have assumed the type of each node is available
as a
covariate. In the stochastic block model, types are latent, and the marginal
distribution of graphs sums over type-conditional distributions. Proposition
1 in the supplementary material~\cite{supp} shows that such summing-over-latents
preserves projectibility. For stochastic block models,
projectibility also follows from~\cite{Lovasz-Szegedy-limits-of-dense-graphs}, Theorem
2.7(ii).} As the distribution
factorizes into a product of $n(n-1)$ terms, each of exactly the same
form, the
log partition function scales exactly with $n(n-1)$, and the
conclusions of Section~\ref{secconsistency-of-mle} imply the strong consistency of the maximum
likelihood estimator.\footnote{An important variant of such models are the
``degree-corrected block models'' of
\cite{Karrer-MEJN-blockmodels-and-community-structure}, where each
node has a
unique parameter, which is its expected degree. It is easily seen that the
range of possible degrees for each new node is the same, no matter what the
configuration of smaller sub-graphs (in which the node does not
appear), as
is the number of configurations giving rise to each degree. The conditions
of Section~\ref{secgrowing-parameter-set} thus hold, and these models are
projective.} This result thus applies to the well-studied $\beta$-model
\cite{Barvinok-Hartigan,Chatterjee-Diaconis-Sly-given-degree-sequence,Rinaldo-Petrovic-Fienberg}.

Typically, however, ERGMs are \emph{not} dyadic independence models. In many
networks, if nodes $i$ and $j$ are both linked to $k$, then $i$ and $j$ are
unusually likely to be directly linked. This will of course happen if
nodes of
the same type are especially likely to be friends (``homophily''
\cite{Birds-of-a-Feather-review}), since then the posterior
probability of $i$
and $j$ being of the same type is elevated. However, it can also be modeled
directly. The direct way to do so is to introduce the number (or
density) of
triangles as a sufficient statistic, but this leads to pathological degeneracy
\cite{Rinaldo-Fienberg-Zhou-network-degeneracy}, and modern specifications
involve a large set of triangle-like motifs
\cite{Snijders-Pattison-Robins-Handcock-new-specification,Wasserman-Robins-dependence-graphs-and-pstar,statnet-special-issue}.
Empirically, when using such specifications, one often finds a nontrivial
coefficient for such ``transitivity'' or ``clustering,'' over and above
homophily~\cite{Goodereau-Kitts-Morris-birds-of-a-feather}. It is
because of
such findings that we ask whether the parameters in these models are
projective.

Sadly, no statistic which counts triangles, or larger motifs, can have
the nice
additive form of dyad counts, no matter how we decompose the network. Take,
for instance, triangles. Any given edge among the first $n$ nodes \emph{could}
be part of a triangle, depending on ties to the next node. Thus to determine
the number of triangles among the first $n+1$ nodes, we need much more
information about the sub-graph of the first $n$ nodes than just the
number of
triangles among them. Indeed, we can go further. The range of possible
increments to the number of triangles changes with the number of existing
triangles. This is quite incompatible with separable increments, so, by
\eqref{thmprojectible-iff-separable}, the parameters cannot be
projective. We remark that the nonprojectibility of Markov graphs
\cite{Frank-Strauss86}, a special instance of ERGMs where the sufficient
statistics count edges, $k$-stars and triangles, was noted in
\cite{Lauritzen-exchangeable-Rasch-matrices}.

Parallel arguments apply to the count of any motif of $k$ nodes, $k >
2$. Any
given edge (or absence of an edge) among the first $n$ nodes could be
part of
such a motif, depending on the edges involving the next $k-2$ nodes. Such
counts are thus not nicely additive. For the same reasons as with triangles,
the range of increments for such statistics is not constant, and nonseparable
increments imply nonprojective family.

While these ERGMs are not projective, some of them may, as a sort of
consolation prize, still satisfy equation \eqref
{eqnasymptotic-scaling}. For
instance, in models where $T$ has two elements, the number of edges and the
(normalized) number of triangles or of 2-stars, the log partition
function is
known to scale like $n(n-1)$ as the number of nodes $n\rightarrow
\infty$, at
least in the parameter regimes where the models behave basically like either
very full or very empty Erd{\H{o}}s--R{\'e}nyi networks
\cite{Park-MEJN-stat-mech-of-networks,Park-MEJN-on-two-star,Park-MEJN-on-clustered,Chatterjee-Dey-on-Steins-method,Chatterjee-Diaconis-ergms,Bhamidi-Bresler-Sly-mixing-time-of-ergms}.
(We suspect, from
\cite{Park-MEJN-stat-mech-of-networks,Xiang-Neville-learning-with-one-network,Chatterjee-Diaconis-ergms},
that similar results apply to many other ERGMs.) Thus, by equation
\eqref{eqnldp-for-suff-stat}, if we fix a large number $n$ of nodes and
generate a graph $X$ from $\mathbb{P}_{\theta,n}$, the probability
that the MLE
$\mle(X)$ will\vadjust{\goodbreak} be more than $\varepsilon$ away from $\theta$ will be
exponentially
small in $n(n-1)$ and $\varepsilon^2$. Since these models are not projective,
however, it is impossible to \emph{improve} parameter estimates by
getting more
data, since parameters for smaller sub-graphs just cannot be
extrapolated to
larger graphs (or vice versa).

We thus have a near-dichotomy for ERGMs. Dyadic independence models have
separable and independent increments to the statistics, and the resulting
family is projective. However, specifications where the sufficient statistics
count larger motifs cannot have separable increments, and
projectibility does
not hold. Such an ERGM may provide a good description of a given social
network on a certain set of nodes, but it cannot be projected to give
predictions on any larger or more global graph from which that one was drawn.
If an ERGM is postulated for the whole network, then inference for its
parameters must explicitly treat the unobserved portions of the network as
missing data (perhaps through an expectation-maximization algorithm),
though of
course there may be considerable uncertainty about just how much data is
missing.

%s6 #&#
\section{Conclusion}

Specifications for exponential families of dependent variables in terms of
joint distributions are surprisingly delicate; the statistics must be chosen
extremely carefully, in order to achieve separable increments. (Conditional
specifications do not have this problem.) This has, perhaps, been
obscured in
the past by the emphasis on using exponential families to model multivariate
but independent cases, as IID models are always projective.

Network models, one of the outstanding applications of exponential families,
suffer from this problem in an acute form. Dyadic independence models are
projective models, but are sociologically extremely implausible, and certainly
do not manage to reproduce the data well. More interesting specifications,
involving clustering terms, never have separable increments. We thus
have an
impasse which it seems can only be resolved by going to a different
family of
specifications. One possibility---which, however, requires more and different
data---is to model the evolution of networks over time~\cite{Snijders-on-longitudinal-network-data}.
In particular,
Hanneke, Fu and Xing~\cite{Hanneke-Fu-Xing-temporal-social-networks}
consider situations where the
distribution of the network at time $t+1$ conditional on the network at time
$t$ follows an exponential family. Even when the statistics in the conditional
specification include (say) changes in the number of triangles, the issues
raised above do not apply.

Roughly speaking, the issue with the nonprojective ERGM
specifications, and
with other nonprojective exponential families, is that the dependency
structure corresponding to the statistics allows interactions between arbitrary
collections of random variables. It is not possible, with those
statistics, to
``screen off'' one part of the assemblage from another by conditioning on
boundary terms. Suppose our larger information set $\largeset$
consists of two
nonoverlapping\vadjust{\goodbreak} and strictly smaller information sets, $\smallset
\subset
\largeset$ and $\alternateset\subset\largeset$, plus the new observation
obtained by looking at both $\smallset$ and $\alternateset$. (For instance,
the latter might be the edges between two disjoint sets of nodes.) Then the
models which work properly are ones where the sufficient statistic for
$\largeset$ partitions into marginal terms from $\smallset$ and
$\alternateset$, plus the interactions strictly between them:
$t_{\largeset}(X_{\largeset}) = t_{\smallset}(X_{\smallset}) +
T_{\alternateset}(X_{\alternateset}) + T_{\largeset\setminus
(\smallset\cup
\alternateset)}(X_{\largeset\setminus(\smallset\cup\alternateset
)})$. In
physical language~\cite{Landau-Lifshitz}, the energy for the whole assemblage
needs to be a sum of two ``volume'' terms for its sub-assemblages, plus a
``surface'' term for their interface. The network models with nonprojective
parameters do not admit such a decomposition; every variable, potentially,
interacts with every other variable.

One might try to give up the exponential family form, while keeping
finite-dimensional sufficient statistics. We suspect that this will not work,
however, since Lauritzen~\cite{Lauritzen-extremal-families-and-suff-stats} showed that
whenever the sufficient statistics form a semi-group, the models must
be either
ordinary exponential families, or certain generalizations thereof with
much the
same properties. We believe that there exists a purely algebraic
characterization of the sufficient statistics compatible with projectibility,
but must leave this for the future.

One reason for the trouble with ERGMs is that every infinite
exchangeable graph
distribution is actually a mixture over projective dyadic-independence
distributions~\cite{Diaconis-Janson-graph-limits,Bickel-Chen-on-modularity},
though not necessarily ones with a finite-dimensional sufficient statistic.
Along any one sequence of sub-graphs from such an infinite graph, in
fact, the
densities of all motifs approach limiting values which pick out a unique
projective dyadic-independence distribution~\cite{Diaconis-Janson-graph-limits}; cf. also
\cite{Lauritzen-extremal-families-and-suff-stats,Lauritzen-exchangeable-Rasch-matrices}.
This suggests that an alternative to parametric inference would be
nonparametric estimation of the limiting dyadic-independence model, by
smoothing the adjacency matrix; this, too, we pursue elsewhere.

\begin{appendix}
%s7 #&#
\section*{Appendix: Proofs}\label{app}

For notation in this section, without loss of generality, fix a generic
pair of
subsets $\smallset\subset\largeset$ and a value of $\theta$. We
will write a
representative point $x_{\largeset} \in\mathcal{X}_{\largeset}$ as
$x_{\largeset} = (x, y)$, with $x \in\mathcal{X}_{\smallset}$ and
$y \in
\mathcal{X}_{\newset}$. Also, we abbreviate $t_{\largeset}(x,y) -
t_{\smallset}(x)$, for $x \in\mathcal{X}_{\smallset}$ and $y \in
\mathcal{X}_{\newset}$ by $t_{\newset}(x,y)$.

%s7.1 #&#
\subsection{\texorpdfstring{Proof of Theorem \protect\ref{thmprojectible-iff-separable}}{Proof of Theorem 1}}

For clarity, we prove the two directions separately. First we show that
projectability implies separable increments.

%pr4 #&#
\begin{proposition}
If the exponential family $\{ \mathbb{P}_\theta\}_{A \in\mathcal
{A}}$ is
projective, then the sufficient statistics $\{ T_A \}_{A \in\mathcal{A}}$
have separable increments, that is, $\smallset\subset\largeset$
implies that
$v_{\newset|\smallset}(\delta,x) = v_{\newset}(\delta)$.
\label{propproj-to-SI}
\end{proposition}

\begin{pf}
By projectibility, for each $\theta$,
%
%e21 #&#
%e22 #&#
%e23 #&#
%e24 #&#
\begin{eqnarray}
p_{\smallset,\theta}(x) & = & \sum_{y \in\mathcal{X}_{\newset
}}{p_{\largeset,\theta}(x,y)
} = \sum_{y \in\mathcal{X}_{\newset}}{ \frac{e^{\langle\theta,t_{\largeset} (x,y) \rangle}}{z_{\largeset}(\theta)} }
\\
& = & \frac{1}{z_{\largeset}(\theta)} \sum_{y \in\mathcal
{X}_{\newset}}{ \exp \bigl\{
\bigl\langle\theta, t_{\largeset} (x,y) - t_{\smallset}(x) \bigr\rangle+
\bigl\langle\theta, t_{\smallset}(x) \bigr\rangle \bigr\}}
\\
& = & \frac{e^{ \langle\theta, t_{\smallset}(x) \rangle
}}{z_{\smallset}(\theta)} \frac{z_{\smallset}(\theta
)}{z_{\largeset}(\theta)} \sum_{y \in\mathcal{X}_{\newset}}
{ \exp \bigl\{ \bigl\langle\theta, t_{\newset} (x,y) \bigr\rangle \bigr\}}
\\
& = & p_{\smallset,\theta}(x) \frac{z_{\smallset}(\theta
)}{z_{\largeset}(\theta)} \sum_{y \in\mathcal{X}_{\newset}}
{ \exp \bigl\{ \bigl\langle\theta, t_{\newset} (x,y) \bigr\rangle \bigr\} },
\end{eqnarray}
which implies that, for all $x \in\mathcal{X}_{\smallset}$,
%
%e25 #&#
\begin{equation}
\sum_{y \in\mathcal{X}_{\newset}}{ \exp \bigl\{ \bigl\langle\theta,
t_{\newset} (x,y) \bigr\rangle \bigr\}} = \frac{z_{\largeset}(\theta
)}{z_{\smallset}(\theta)}.
\label{eqconstant-conditional-partition-function}
\end{equation}
Re-writing the left-hand side of equation
\eqref{eqconstant-conditional-partition-function} as a sum over the set
$\Delta(x)$ of values which the increment $t_{\newset} (x,y)$ to the
sufficient
statistic might take yields
%
%e26 #&#
\begin{equation}
\sum_{\delta\in\Delta(x)}{v_{\newset|\smallset}(\delta,x)\exp {\langle
\theta, \delta\rangle}} = \frac{z_{\largeset}(\theta
)}{z_{\smallset}(\theta)}, \label{eqnconditional-partition-function-as-sum-over-increment-values}
\end{equation}
where the joint volume factor is defined in \eqref
{eqjointvolfactor}. Since the right-hand side of equation
\eqref
{eqnconditional-partition-function-as-sum-over-increment-values} is the
same for all $x$, so must the left-hand side.

Observe that this left-hand side is the Laplace transform of the function
$v_{\newset|\smallset}(\cdot,x)$. The latter is a nonnegative
function which
defines a measure on $\mathbb{R}^d$, whose support is $\Delta(x)$. Hence,
%
%e27 #&#
\begin{equation}
\sum_{\delta\in\Delta(x)}{\frac{v_{\newset|\smallset}(\delta,x)}{\sum_{\delta^{\prime}\in\Delta(x)}{v_{\newset|\smallset
}(\delta^{\prime},x)}}\exp{\langle\theta,
\delta\rangle}} \label{eqnlaplace-transform-is-constant}
\end{equation}
is the Laplace transform of a discrete probability measure in $\mathbb{R}^d$.
But the denominator in the inner sum is just $|\mathcal{X}_{\newset
}|$, no
matter what $x$ might be.\footnote{This can be seen either from
recalling that
exponential families have full support, or from defining $T_{\largeset
}$ as a
total and not a partial function on $\mathcal{X}_{\largeset}$.} So
we have
that for any $x, x^{\prime} \in\mathcal{X}_{\smallset}$, and all
$\theta\in
\Theta$,
%
%e28 #&#
\begin{equation}
\sum_{\delta\in\Delta(x)}{\frac{v_{\newset|\smallset}(\delta,x)}{|\mathcal{X}_{\newset}|}\exp{\langle\theta,
\delta\rangle}} = \sum_{\delta\in\Delta(x^{\prime})}{\frac{v_{\newset|\smallset
}(\delta,x^{\prime})}{|\mathcal{X}_{\newset}|}\exp{
\langle\theta, \delta\rangle}}. \label{eqntwo-laplace-transforms}
\end{equation}
Since both sides of equation \eqref{eqntwo-laplace-transforms} are Laplace
transforms of probability measures on a common space, and the equality
holds on
all of $\Theta$, which contains an open set, we may conclude that the two
measures are equal~\cite{Barndorff-Nielsen-exponential-families}, Theorem~7.3.
This means that they have the same support, $\Delta(x) = \Delta
(x^{\prime}) =
\Delta$, and that they have the same density with respect to counting measure
on~$\Delta$. As they also have the same normalizing factor (viz.,
$|\mathcal{X}_{\newset}|$), we get that $v_{\newset|\smallset
}(\delta,x) =
v_{\newset|\smallset}(\delta,x^{\prime}) = v_{\newset}(\delta)$.
Since the
points $x$ and $x'$ are arbitrary, this last property is precisely having
separable increments.
\end{pf}

Next, we prove the reverse direction, namely that separable increments imply
projectibility. This is clearer with some preliminary lemmas.

%le1 #&#
\begin{lemma}\label{lemmaseparable-increments-implies-volume-factor-separation}
If the sufficient statistics have separable increments, then the joint volume
factors factorize, that is,
%
%e29 #&#
\begin{equation}
v_{\smallset,\newset}(t,\delta) = v_{\smallset}(t) v_{\newset
}(\delta)
\end{equation}
for all $A \subset B$, $t$ and $\delta$.
\end{lemma}
\begin{pf}
By definition,
%
%e30 #&#
\begin{equation}
v_{\smallset,\newset}(t,\delta) = \sum_{\{ x\in\mathcal
{X}_{\smallset} \colon t_{\smallset}(x) = t
\}}{v_{\newset|\smallset}(
\delta,x)}. \label{eqnjoint-volume-factor-is-sum-of-conditional-volume-factors}
\end{equation}
When the statistic has separable increments, $v_{\newset|\smallset
}(\delta,x)
= v_{\newset}(\delta)$, so
%
%e31 #&#
\begin{equation}
v_{\smallset,\newset}(t,\delta) = \sum_{\{ x \colon t_{\smallset
}(x) = t
\}}{v_{\newset}(
\delta)} = v_{\smallset}(t) v_{\newset}(\delta),
\end{equation}
proving the claim.
\end{pf}

%le2 #&#
\begin{lemma}\label{lemmaseparable-volume-factor-implies-projectible-increments}
If the joint volume factor factorizes, then the sufficient statistics has
independent increments, and the distribution of the sufficient static is
projective.
\end{lemma}
\begin{pf} Without loss of generality, fix a value $t$ for
$T_{\smallset}$
and $\delta$ for $T_{\newset}$. By the law of total probability and the
definition of the volume factor,
%
%e32 #&#
\begin{equation}
\mathbb{P}_{\theta,\largeset}(T_{\smallset}=t,T_{\newset}=\delta) =
v_{\smallset,\newset}(t,\delta) \frac{e^{\langle\theta, t \rangle}
e^{\langle\theta, \delta\rangle}}{z_\largeset(\theta)}.
\end{equation}
If the volume factor factorizes, so that $v_{\smallset,\newset
}(t,\delta) =
v_{\smallset}(t) v_{\newset}(\delta)$, then we obtain
%
%e33 #&#
\begin{equation}\quad
\mathbb{P}_{\theta,\largeset}(T_{\smallset}=t,T_{\newset}=\delta) = \biggl[
\frac{1}{z_{A}(\theta)} v_{\smallset}(t) e^{\langle\theta, t
\rangle} \biggr] \biggl[
\frac{z_{A}(\theta)}{z_{B}(\theta)} v_{\newset}(\delta) e^{\langle\theta,\delta\rangle} \biggr].
\end{equation}
It then follows that
%
%e34 #&#
\begin{equation}
\mathbb{P}_{\theta,\largeset}(T_{\smallset}=t,T_{\newset}=\delta) =
\mathbb{P}_{\theta,\largeset}(T_{\smallset}=t)\mathbb{P}_{\theta,\largeset}(T_{\newset}=
\delta)\qquad \forall\theta,
\end{equation}
and thus that $T$ has independent increments. To establish the
projectibility of the distribution of $T$, sum over $\delta$
\begin{eqnarray*}
\mathbb{P}_{\theta,\largeset}(T_{\smallset}=t) & = & \sum
_{\delta
}{\mathbb{P}_{\theta,\largeset}(T_{\smallset}=t,T_{\newset}=
\delta )}
\\
&=& \frac{ v_{\smallset}(t) e^{\langle\theta, t \rangle
}}{z_{\largeset}(\theta)}\sum_{\delta}{v_{\newset}(
\delta) e^{\langle\theta,\delta\rangle}}
\\
& = & \frac{v_{\smallset}(t) e^{\langle\theta, t \rangle
}}{z_{\largeset}(\theta)} z_{\newset}(\theta).
\end{eqnarray*}
Since $\mathbb{P}_{\smallset,\theta}(T_{\smallset}=t) =
v_{\smallset}(t)e^{\langle\theta, t \rangle}/z_{\smallset}(\theta
)$, and
both distributions must sum to 1 over~$t$, we can conclude that
$z_{\smallset}(\theta) = z_{\largeset}(\theta)/z_{\newset}(\delta
)$, and
hence that the distribution of the sufficient statistic is projective.
\end{pf}

%le3 #&#
\begin{lemma}%
\label{lemmamixed-joint-prob-of-obs-and-stat}
If the sufficient statistics of an exponential family have separable
increments, then
%
%e35 #&#
\begin{equation}\quad
\mathbb{P}_{\largeset,\theta}(X_{\smallset}=x,T_{\newset}=\delta) =
\frac{1}{v_{\smallset}(t_{\smallset}(x))}\mathbb{P}_{\largeset,\theta}\bigl(T_{\smallset}=t_{\smallset}(x),T_{\newset}=
\delta\bigr).
\end{equation}
\end{lemma}
\begin{pf} Abbreviate $t_{\smallset}(x)$ by $t$. By the law of total
probability,
%
%e36 #&#
\begin{equation}
\mathbb{P}_{\largeset,\theta}(T_{\smallset}=t,T_{\newset}=\delta) = \sum
_{(x,y)\dvtx t_{\smallset}(x) = t, t_{\newset}(x,y)
=\delta}{p_{\largeset,\theta}(x,y)}.
\end{equation}
Since $T_{\largeset}$ is sufficient, and $t_{\largeset}(x,y) =
t+\delta$ for
all $(x,y)$ in the sum,
%
%e37 #&#
\begin{equation}
\mathbb{P}_{\largeset,\theta}(T_{\smallset}=t,T_{\newset}=\delta) =
v_{\smallset,\newset}(t,\delta) e^{\langle\theta, t+\delta\rangle
} / z_{\largeset}(\theta).
\end{equation}
By parallel reasoning,
%
%e38 #&#
\begin{equation}
\mathbb{P}_{\largeset,\theta}(X_{\smallset}=x,T_{\newset}=\delta) =
v_{\newset|\smallset}(\delta,x) e^{\langle\theta, t+\delta\rangle
} / z_{\largeset}(\theta).
\end{equation}
Therefore,
%
%e39 #&#
\begin{equation}
\mathbb{P}_{\largeset,\theta}(X_{\smallset}=x,T_{\newset}=\delta) =
v_{\newset|\smallset}(\delta,x) \frac{\mathbb{P}_{\largeset,\theta}(T_{\smallset}=t,T_{\newset
}=\delta)}{v_{\smallset,\newset}(t,\delta)}.
\end{equation}
If the statistic has separable increments, then
$v_{\smallset,\newset}(t,\delta) = v_{\smallset}(t)v_{\newset
}(\delta) =
v_{\smallset}(t)v_{\newset|\smallset}(\delta,x)$, and the
conclusion follows.
\end{pf}

\begin{remark*}
The lemma does \emph{not} follow merely from the joint volume
factor separating, $v_{\smallset,\newset}(t,\delta) =
v_{\smallset}(t)v_{\newset}(\delta)$. The conditional volume factor
must also
be constant in $x$.
\end{remark*}

%pr5 #&#
\begin{proposition}
If the sufficient statistic of an exponential family has separable
increments, then the family is projective.
\label{propSI-to-proj}
\end{proposition}

\begin{pf} We calculate the marginal probability of $X_{\smallset}$ in
$\mathbb{P}_{\theta,\largeset}$, by integrating out the increment to the
sufficient statistic. (The set of possible increments, $\Delta$, is
the same
for all $x$, by separability.) Once again, we abbreviate $t_{\smallset}(x)$
by $t$:
\begin{eqnarray*}
\mathbb{P}_{\largeset,\theta}(X_{\smallset}=x) & = & \sum
_{\delta
\in\Delta}{\mathbb{P}_{\largeset,\theta}(X_{\smallset}=x,
T_{\newset}=\delta)}
\\
& = & \frac{1}{v_{\smallset}(t)}\sum_{\delta\in\Delta}{\mathbb
{P}_{\largeset,\theta}(T_{\smallset}=t, T_{\newset} = \delta)}
\\
& = & \frac{1}{v_{\smallset}(t)}\sum_{\delta\in\Delta}{\mathbb
{P}_{\largeset,\theta}(T_{\smallset}=t)\mathbb{P}_{\largeset,\theta}(T_{\newset}
= \delta|T_{\smallset}=t)}
\\
& = & \frac{\mathbb{P}_{\largeset,\theta}(T_{\smallset
}=t)}{v_{\smallset}(t)}
\\
& = & \frac{\mathbb{P}_{\smallset,\theta}(T_{\smallset
}=t)}{v_{\smallset}(t)}
\\
& = & \mathbb{P}_{\smallset,\theta}(X_{\smallset}=x).
\end{eqnarray*}
These steps use, in succession: Lemma
\ref{lemmamixed-joint-prob-of-obs-and-stat}; the fact that conditional
probabilities sum to 1; the projectibility of the sufficient statistics (via
Lemmas \ref
{lemmaseparable-increments-implies-volume-factor-separation} and
\ref{lemmaseparable-volume-factor-implies-projectible-increments});
and the
definition of $v_{\smallset}(t)$.
\end{pf}

%s7.2 #&#
\subsection{Other proofs}
\mbox{}
\begin{pf*}{Proof of Proposition~\ref{propproj-to-II}} By Proposition
\ref{propproj-to-SI}, a projective family has separable increments,
and by
Lemma~\ref{lemmaseparable-volume-factor-implies-projectible-increments},
separable increments implies independent increments.
\end{pf*}

\begin{pf*}{Proof of Proposition
\ref{propincrement-indep-of-old-observable}}
By Proposition~\ref{propproj-to-SI}, every projective exponential
family has
separable increments. By Lemma~\ref{lemmamixed-joint-prob-of-obs-and-stat},
in an exponential family with separable increments,
%
%e40 #&#
\begin{equation}\quad
\mathbb{P}_{\largeset,\theta}(X_{\smallset}=x,T_{\newset}=\delta) =
\frac{1}{v_{\smallset}(t_{\smallset}(x))}\mathbb{P}_{\largeset,\theta}\bigl(T_{\smallset}=t_{\smallset}(x),T_{\newset}=
\delta\bigr).
\end{equation}
Therefore, using projectibility,
%
%e41 #&#
\begin{equation}\quad
\mathbb{P}_{\largeset,\theta}(T_{\newset}=\delta|X_{\smallset}=x) =
\frac{\mathbb{P}_{\largeset,\theta}(T_{\smallset}=t_{\smallset
}(x),T_{\newset}=\delta)/v_{\smallset}(t_{\smallset
}(x))}{p_{\smallset,\theta}(x)}.
\end{equation}
By the definition of $v_{\smallset}(\cdot)$, $p_{\smallset,\theta
}(x) =
\mathbb{P}_{\smallset,\theta}(T_{\smallset}=t_{\smallset
}(x))/v_{\smallset}(t_{\smallset}(x))$,
so
%
%e42 #&#
\begin{equation}\quad
\mathbb{P}_{\largeset,\theta}(T_{\newset}=\delta|X_{\smallset}=x) =
\frac{\mathbb{P}_{\largeset,\theta}(T_{\smallset}=t_{\smallset
}(x),T_{\newset}=\delta)}{\mathbb{P}_{\smallset,\theta
}(T_{\smallset}=t_{\smallset}(x))}.
\end{equation}
But, by Lemma
\ref{lemmaseparable-volume-factor-implies-projectible-increments}, the
sufficient statistics have a projective distribution with independent
increments, implying
%
%e43 #&#
\begin{equation}\quad
\mathbb{P}_{\largeset,\theta}\bigl(T_{\smallset}=t_{\smallset
}(x),T_{\newset}=
\delta\bigr) = \mathbb{P}_{\smallset,\theta}\bigl(T_{\smallset}=t_{\smallset
}(x)
\bigr)\mathbb{P}_{\largeset,\theta}(T_{\newset}=\delta).
\end{equation}
Therefore,
%
%e44 #&#
\begin{equation}
\mathbb{P}_{\largeset,\theta}(T_{\newset}=\delta|X_{\smallset}=x) =
\mathbb{P}_{\largeset,\theta}(T_{\newset}=\delta)
\end{equation}
and so $T_{\newset} \indep X_{\smallset}$.
\end{pf*}

\begin{pf*}{Proof of Proposition
\ref{propII-to-volume-factor-separation}}
Below we prove that if the suffiicient statistics of an exponential family
have independent increments, then the volume factor separates, and the
distribution of the statistic is projective.

Since $T_{\largeset}$ is a sufficient statistic, by the Neyman factorization
theorem (\cite{Schervish-theory-of-stats}, Theorem 2.21, page 89),
%
%e45 #&#
\begin{equation}
\mathbb{P}_{\largeset,\theta}(X_{\smallset}=x,X_{\newset}=y) =
g_{\largeset}\bigl(\theta,t_{\smallset}(x)+t_{\newset}(x,y)\bigr)
h(x,y).
\end{equation}
In light of equation \eqref{eqnexponential-family-density}, we may take
$h(x,y) = 1$. Abbreviating $t_{\smallset}(x)$ by $t$ and $t_{\newset}(x,y)$
by $\delta$, it follows that
%
%e46 #&#
\begin{equation}
\mathbb{P}_{\largeset,\theta}(T_{\smallset}=t,T_{\newset}=\delta) =
v_{\smallset,\newset}(t,\delta)g_{\largeset}(\theta,t+\delta) \label{eqnequiprobability}.
\end{equation}

By independent increments, however,
%
%e47 #&#
\begin{equation}
\mathbb{P}_{\largeset,\theta}(T_{\smallset}=t,T_{\newset}=\delta) =
\mathbb{P}_{\largeset,\theta}(T_{\smallset}=t)\mathbb {P}_{\largeset,\theta}(T_{\newset}=
\delta),
\end{equation}
whence it follows that, for some functions $g_{\newset}, k_{\smallset},
k_{\newset}$,
%
%e48 #&#
\begin{equation}
g_{\largeset}(\theta,t+\delta) = g_{\smallset}(\theta,t) g_{\newset}(
\theta,\delta)
\end{equation}
and
%
%e49 #&#
\begin{equation}
v_{\smallset,\newset}(t,\delta) = k_{\smallset}(t)k_{\newset
}(\delta)
\end{equation}
and
%
%e50 #&#
\begin{equation}
\mathbb{P}_{\largeset,\theta}(T_{\smallset}=t,T_{\newset}=\delta) =
k_{\smallset}(t)k_{\newset}(\delta)g_{\smallset}(\theta,t)
g_{\newset}(\theta,\delta).
\end{equation}
To proceed, we must identify the new $g$ and $k$ functions. To this end,
recalling that $v_{\smallset}(t)$ is the number of $x_{\smallset}$
configurations such that $t_{\smallset}(x_{\smallset}) = t$, we have
%
%e51 #&#
\begin{equation}
\sum_{\delta}{v_{\smallset,\newset}(t,\delta)} =
v_{\smallset}(t)|\mathcal{X}_{\newset}|
\end{equation}
and, at the same time,
%
%e52 #&#
\begin{equation}
\sum_{\delta}{v_{\smallset,\newset}(t,\delta)} =
k_{\smallset}(t)\sum_{\delta}{k_{\newset}(
\delta)}.
\end{equation}
Clearly, then, $k_{\smallset}(t) = c_1 v_{\smallset}(t)$ while
$\sum_{\delta}{k_{\newset}(\delta)} = c_2|\mathcal{X}_{\newset
}|$. Since
%
%e53 #&#
\begin{equation}
\sum_{t}{\sum_{\delta}{v_{\smallset,\newset}(t,
\delta)}} = |\mathcal{X}_{\smallset}||\mathcal{X}_{\newset}|
\end{equation}
and $\sum_{t}{v_{\smallset}(t)} =|\mathcal{X}_{\smallset}|$, we
need $c_1 c_2
= 1$, and may take $c_1 = c_2 =1 $ for simplicity. This allows us to write
%
%e54 #&#
\begin{equation}
v_{\smallset,\newset}(t,\delta) = v_{\smallset}(t)v_{\newset
}(\delta),
\end{equation}
which is exactly the assertion that the volume factor separates.

Turning to the $g$ functions, we sum over $\delta$ again to obtain the
marginal distribution of $T_{\smallset}$,
\begin{eqnarray*}
\mathbb{P}_{\largeset,\theta}(T_{\smallset} = t) & = & \sum
_{\delta}{\mathbb{P}_{\largeset,\theta}(T_{\smallset}=t,T_{\newset
}=
\delta)}
\\
& = & \sum_{\delta}{v_{\smallset}(t)
g_{\smallset}(\theta,t) v_{\newset}(\delta) g_{\newset}(\theta,
\delta)}
\\
& = & v_{\smallset}(t) g_{\smallset}(\theta,t)\sum
_{\delta
}{v_{\newset}(\delta)g_{\newset}(\theta,
\delta)}.
\end{eqnarray*}
Now, we finally we use the exponential-family form. Specifically, we know
that
%
%e55 #&#
\begin{equation}
g_{\largeset}(\theta,t+\delta) = \frac{e^{\langle\theta, t\rangle
}e^{\langle
\theta, \delta\rangle}}{z_{\largeset}(\theta)},
\end{equation}
so that $g_{\smallset}(\theta,t) \propto e^{\langle\theta, t\rangle}$,
$g_{\newset}(\theta,\delta) \propto e^{\langle\theta, \delta
\rangle}$.
Therefore,
%
%e56 #&#
\begin{equation}
\mathbb{P}_{\largeset,\theta}(T_{\smallset} = t) \propto v_{\smallset}(t)e^{\langle\theta, t\rangle}
\propto \mathbb{P}_{\smallset,\theta}(T_{\smallset} = t),
\end{equation}
and normalization now forces
%
%e57 #&#
\begin{equation}
\mathbb{P}_{\largeset,\theta}(T_{\smallset} = t) = \mathbb{P}_{\smallset,\theta}(T_{\smallset}
= t)
\end{equation}
as desired.
\end{pf*}

\begin{pf*}{Proof of Theorem~\ref{thmpredictive-sufficiency}}
The conditional density of $X_{\newset}$ given $X_{\smallset}$ is
just the
ratio of joint to marginal densities (both with the same $\theta$, by
projectibility),
%
%e58 #&#
%e59 #&#
\begin{eqnarray}
p_{\largeset|\smallset,\theta}(y|x) & = & \frac{p_{\largeset,\theta}(x,y)}{p_{\smallset,\theta}(x)} = \frac{e^{\langle\theta, t_{\largeset}(x,y)\rangle} /
z_{\largeset}(\theta)}{e^{\langle\theta, t_\smallset(x)\rangle}/
z_\smallset(\theta)}
\\
& = & \frac{e^{\langle\theta, t_{\newset}(x,y) \rangle
}}{z_{\largeset}(\theta)/z_\smallset(\theta)},
\end{eqnarray}
which is an exponential family with parameter $\theta$, sufficient statistic
$T_{\newset}$, and partition function $z_{\newset|\smallset}(\theta
) \equiv
z_{\largeset}(\theta)/z_\smallset(\theta)$.
\end{pf*}

\begin{pf*}{Proof of Theorem~\ref{thmldp}}
Under equation \eqref{eqnasymptotic-scaling}, the cumulant generating
function also scales asymptotically,\vadjust{\goodbreak}
$\kappa_{\smallset,\theta}(\phi)/r_{|\smallset|} \rightarrow
a(\theta+\phi) -
a(\theta)$. Since $a$ is differentiable, the G{\"a}rtner--Ellis
theorem of
large deviations theory~\cite{den-Hollander-large-deviations}, Chapter V, implies
that $T_{\smallset}/r_{|\smallset|}$ obeys a large deviations
principle with
rate $r_{|\smallset|}$, and rate function given by equation
\eqref{eqnrate-function}, which is to say, equations~\eqref{eqnldp-lower}
and \eqref{eqnldp-upper}.
\end{pf*}

%le4 #&#
\begin{lemma}\label{corconditional-moment-generating-function}
The moment generating function of $T_{\newset}$ is
%
%e60 #&#
\begin{equation}
\frac{z_{\largeset}(\theta+\phi) z_\smallset(\theta
)}{z_{\largeset}(\theta) z_\smallset(\theta+\phi)} = \frac
{M_{\theta,\largeset}(\phi)}{M_{\theta,\smallset}(\phi)} \label{eqnconditional-moment-generating-function}.
\end{equation}
\end{lemma}

\begin{pf}
From the proof of Theorem~\ref{thmpredictive-sufficiency},
$X_{\newset}|X_{\smallset}$ has an exponential family distribution with
sufficient statistic $T_{\newset}$. Thus we may use equation
\eqref{eqnmgf-of-exponential-family} to find the moment generating function
of $T_{\newset}$ conditional on $X_{\smallset}$,
%
%e61 #&#
%e62 #&#
%e63 #&#
\begin{eqnarray}
M_{\theta, \newset|\smallset}(\phi) & = & \frac{z_{\newset
|\smallset}(\theta+\phi)}{z_{\newset|\smallset}(\theta)}
\\
& = & \frac{z_{\largeset}(\theta+\phi)/z_\smallset(\theta+\phi
)}{z_{\largeset}(\theta)/z_\smallset(\theta)}
\\
& = & \frac{z_{\largeset}(\theta+\phi) z_\smallset(\theta
)}{z_{\largeset}(\theta) z_\smallset(\theta+\phi)} = \frac{M_{\theta,\largeset}(\phi)}{M_{\theta,\smallset}(\phi
)}.
\end{eqnarray}
Since, however, $T_{\newset} \indep X_{\smallset}$ (Proposition
\ref{propincrement-indep-of-old-observable}), equation
\eqref{eqnconditional-moment-generating-function} must also give the
unconditional moment generating function.
\end{pf}
\end{appendix}

\section*{Acknowledgments}
We thank Luis Carvalho, Aaron Clauset, Mark
Handcock, Steve Hanneke, Brain Karrer, Sergey Kirshner, Steffen Lauritzen,
David Lazer, John Miller, Martina Morris, Jennifer Neville, Mark
Newman, Peter
Orbanz, Andrew Thomas and Chris Wiggins, for valuable conversations; an
anonymous referee of an earlier version for pointing out a gap in a
proof; and
audiences at the Boston University probability and statistics seminar, and
Columbia University's applied math seminar.

\begin{supplement}[id=suppA]
\stitle{Non-uniform base measures and conditional projectibility\\}
\slink[doi]{10.1214/12-AOS1044SUPP} %[doi,] - jei reikia suskaldyti doi
\sdatatype{.pdf}
\sfilename{aos1044\_supp.pdf}
\sdescription{In the supplementary material we consider the case of
nonuniform base measures and also study a more general form of
conditional projectibility, which implies, in particular, that
stochastic block models are projective.}
\end{supplement}

% imsref loaded by akundreckaite, 2013-03-20 08:38:02
% imsref loaded by akundreckaite, 2013-03-20 11:28:20

\printaddresses


\begin{thebibliography}{69}
% BibTex style file: ims.bst, 2013-01-28
% Default style options (sort=0,type=number).
% Used options (sort=1,type=number).

\bibitem{Achlioptas-et-al-bias-of-traceroute}
\begin{binproceedings}[author]
\bauthor{\bsnm{Achlioptas},~\bfnm{Dimitris}\binits{D.}},
  \bauthor{\bsnm{Clauset},~\bfnm{Aaron}\binits{A.}},
  \bauthor{\bsnm{Kempe},~\bfnm{David}\binits{D.}} \AND
  \bauthor{\bsnm{Moore},~\bfnm{Cristopher}\binits{C.}}
(\byear{2005}).
\btitle{On the bias of traceroute sampling (or: {W}hy almost every network
  looks like it has a power law)}.
In \bbooktitle{Proceedings of the 37th ACM Symposium on Theory of Computing}.
\bptok{imsref}%
\end{binproceedings}
\endbibitem

\bibitem{Ackland-ONeil-online-collective-identity}
\begin{barticle}[author]
\bauthor{\bsnm{Ackland},~\bfnm{Robert}\binits{R.}} \AND
  \bauthor{\bsnm{O'Neil},~\bfnm{Mathieu}\binits{M.}}
(\byear{2011}).
\btitle{Online collective identity: The case of the environmental movement}.
\bjournal{Social Networks}
\bvolume{33}
\bpages{177--190}.
\bptok{imsref}%
\end{barticle}
\endbibitem

\bibitem{Ahmed-Neville-Kompella-network-sampling}
\begin{binproceedings}[author]
\bauthor{\bsnm{Ahmed},~\bfnm{Nesreen~K.}\binits{N.~K.}},
  \bauthor{\bsnm{Neville},~\bfnm{Jennifer}\binits{J.}} \AND
  \bauthor{\bsnm{Kompella},~\bfnm{Ramana}\binits{R.}}
(\byear{2010}).
\btitle{Reconsidering the foundations of network sampling}.
In \bbooktitle{Proceedings of the 2nd Workshop on Information in Networks [WIN
  2010]}
(\beditor{\bfnm{Sinan}\binits{S.}~\bsnm{Aral}},
  \beditor{\bfnm{Foster}\binits{F.}~\bsnm{Provost}} \AND
  \beditor{\bfnm{Arun}\binits{A.}~\bsnm{Sundararajan}}, eds.).
\bptok{imsref}%
\end{binproceedings}
\endbibitem

\bibitem{pstar-primer}
\begin{barticle}[author]
\bauthor{\bsnm{Anderson},~\bfnm{Carolyn~J.}\binits{C.~J.}},
  \bauthor{\bsnm{Wasserman},~\bfnm{Stanley}\binits{S.}} \AND
  \bauthor{\bsnm{Crouch},~\bfnm{Bradley}\binits{B.}}
(\byear{1999}).
\btitle{A $p^*$ primer: Logit models for social networks}.
\bjournal{Social Networks}
\bvolume{21}
\bpages{37--66}.
\bptok{imsref}%
\end{barticle}
\endbibitem

\bibitem{Bahadur-limit-theorems}
\begin{bbook}[author]
\bauthor{\bsnm{Bahadur},~\bfnm{R.~R.}\binits{R.~R.}}
(\byear{1971}).
\btitle{Some Limit Theorems in Statistics}.
\bpublisher{SIAM}, \blocation{Philadelphia}.
\bptok{imsref}%
\end{bbook}
\endbibitem

\bibitem{Barndorff-Nielsen-exponential-families}
\begin{bbook}[mr]
\bauthor{\bsnm{Barndorff-Nielsen},~\bfnm{Ole}\binits{O.}}
(\byear{1978}).
\btitle{Information and Exponential Families in Statistical Theory}.
\bpublisher{Wiley}, \blocation{Chichester}.
\bid{mr={0489333}}
\bptok{imsref}%
\end{bbook}
\endbibitem

\bibitem{Barvinok-Hartigan}
\begin{bmisc}[author]
\bauthor{\bsnm{Barvinok},~\bfnm{A.}\binits{A.}} \AND
  \bauthor{\bsnm{Hartigan},~\bfnm{J.~A.}\binits{J.~A.}}
(\byear{2010}).
\bhowpublished{The number of graphs and a random graph with a given degree
  sequence. Available at arXiv:\arxivurl{1003.0356}.}
\bptok{imsref}%
\end{bmisc}
\endbibitem

\bibitem{Besag-candidates-formula}
\begin{barticle}[mr]
\bauthor{\bsnm{Besag},~\bfnm{Julian}\binits{J.}}
(\byear{1989}).
\btitle{A candidate's formula: A curious result in {B}ayesian prediction}.
\bjournal{Biometrika}
\bvolume{76}
\bpages{183}.
\bid{doi={10.1093/biomet/76.1.183}, issn={0006-3444}, mr={0991437}}
\bptok{imsref}%
\end{barticle}
\endbibitem

\bibitem{Bhamidi-Bresler-Sly-mixing-time-of-ergms}
\begin{barticle}[mr]
\bauthor{\bsnm{Bhamidi},~\bfnm{Shankar}\binits{S.}},
  \bauthor{\bsnm{Bresler},~\bfnm{Guy}\binits{G.}} \AND
  \bauthor{\bsnm{Sly},~\bfnm{Allan}\binits{A.}}
(\byear{2011}).
\btitle{Mixing time of exponential random graphs}.
\bjournal{Ann. Appl. Probab.}
\bvolume{21}
\bpages{2146--2170}.
\bid{issn={1050-5164}, mr={2895412}}
\bptok{imsref}%
\end{barticle}
\endbibitem

\bibitem{Bickel-Chen-on-modularity}
\begin{barticle}[author]
\bauthor{\bsnm{Bickel},~\bfnm{Peter~J.}\binits{P.~J.}} \AND
  \bauthor{\bsnm{Chen},~\bfnm{Aiyou}\binits{A.}}
(\byear{2009}).
\btitle{A nonparametric view of network models and {Newman}--{Girvan} and other
  modularities}.
\bjournal{Proc. Natl. Acad. Sci. USA}
\bvolume{106}
\bpages{21068--21073}.
\bptok{imsref}%
\end{barticle}
\endbibitem

\bibitem{Brown-on-exponential-families}
\begin{bbook}[mr]
\bauthor{\bsnm{Brown},~\bfnm{Lawrence~D.}\binits{L.~D.}}
(\byear{1986}).
\btitle{Fundamentals of Statistical Exponential Families with Applications in
  Statistical Decision Theory}.
\bseries{Institute of Mathematical Statistics Lecture Notes---Monograph Series}
\bvolume{9}.
\bpublisher{IMS}, \blocation{Hayward, CA}.
\bid{mr={0882001}}
\bptok{imsref}%
\end{bbook}
\endbibitem

\bibitem{Butler-predictive-likelihood}
\begin{barticle}[mr]
\bauthor{\bsnm{Butler},~\bfnm{Ronald~W.}\binits{R.~W.}}
(\byear{1986}).
\btitle{Predictive likelihood inference with applications}.
\bjournal{J. Roy. Statist. Soc. Ser. B}
\bvolume{48}
\bpages{1--38}.
\bid{issn={0035-9246}, mr={0848048}}
\bptnote{check related}%
\bptok{imsref}%
\end{barticle}
\endbibitem

%  \beditor{\bsnm{Scott},~\bfnm{John}\binits{J.}} \AND
%  \beditor{\bsnm{Wasserman},~\bfnm{Stanley}\binits{S.}}, eds.
%(\byear{2005}).

\bibitem{Chatterjee-Dey-on-Steins-method}
\begin{barticle}[mr]
\bauthor{\bsnm{Chatterjee},~\bfnm{Sourav}\binits{S.}} \AND
  \bauthor{\bsnm{Dey},~\bfnm{Partha~S.}\binits{P.~S.}}
(\byear{2010}).
\btitle{Applications of {S}tein's method for concentration inequalities}.
\bjournal{Ann. Probab.}
\bvolume{38}
\bpages{2443--2485}.
\bid{doi={10.1214/10-AOP542}, issn={0091-1798}, mr={2683635}}
\bptok{imsref}%
\end{barticle}
\endbibitem

\bibitem{Chatterjee-Diaconis-ergms}
\begin{bmisc}[author]
\bauthor{\bsnm{Chatterjee},~\bfnm{Sourav}\binits{S.}} \AND
  \bauthor{\bsnm{Diaconis},~\bfnm{Persi}\binits{P.}}
(\byear{2011}).
\bhowpublished{Estimating and understanding exponential random graph models.
  Available at arXiv:\arxivurl{1102.2650}.}
\bptok{imsref}%
\end{bmisc}
\endbibitem

\bibitem{Chatterjee-Diaconis-Sly-given-degree-sequence}
\begin{barticle}[author]
\bauthor{\bsnm{Chatterjee},~\bfnm{Sourav}\binits{S.}},
  \bauthor{\bsnm{Diaconis},~\bfnm{Persi}\binits{P.}} \AND
  \bauthor{\bsnm{Sly},~\bfnm{Allan}\binits{A.}}
(\byear{2011}).
\btitle{Random graphs with a given degree sequence}.
\bjournal{Ann. Appl. Probab.}
\bvolume{21}
\bpages{1400--1435}.
\bid{mr={2857452}}
\bptok{imsref}%
\end{barticle}
\endbibitem

\bibitem{Daraganova-et-al-networks-and-geography}
\begin{barticle}[author]
\bauthor{\bsnm{Daraganova},~\bfnm{Galina}\binits{G.}},
  \bauthor{\bsnm{Pattison},~\bfnm{Pip}\binits{P.}},
  \bauthor{\bsnm{Koskinen},~\bfnm{Johan}\binits{J.}},
  \bauthor{\bsnm{Mitchell},~\bfnm{Bill}\binits{B.}},
  \bauthor{\bsnm{Bill},~\bfnm{Anthea}\binits{A.}},
  \bauthor{\bsnm{Watts},~\bfnm{Martin}\binits{M.}} \AND
  \bauthor{\bsnm{Baum},~\bfnm{Scott}\binits{S.}}
(\byear{2012}).
\btitle{Networks and geography: Modelling community network structure as the
  outcome of both spatial and network processes}.
\bjournal{Social Networks}
\bvolume{34}
\bpages{6--17}.
\bptok{imsref}%
\end{barticle}
\endbibitem

\bibitem{de-la-Haye-et-al-on-obesity-and-friendship-networks}
\begin{barticle}[author]
\bauthor{\bparticle{de~la} \bsnm{Haye},~\bfnm{Kayla}\binits{K.}},
  \bauthor{\bsnm{Robins},~\bfnm{Garry}\binits{G.}},
  \bauthor{\bsnm{Mohr},~\bfnm{Philip}\binits{P.}} \AND
  \bauthor{\bsnm{Wilson},~\bfnm{Carlene}\binits{C.}}
(\byear{2010}).
\btitle{Obesity-related behaviors in adolescent friendship networks}.
\bjournal{Social Networks}
\bvolume{32}
\bpages{161--167}.
\bptok{imsref}%
\end{barticle}
\endbibitem

\bibitem{den-Hollander-large-deviations}
\begin{bbook}[mr]
\bauthor{\bparticle{den} \bsnm{Hollander},~\bfnm{Frank}\binits{F.}}
(\byear{2000}).
\btitle{Large Deviations}.
\bseries{Fields Institute Monographs}
\bvolume{14}.
\bpublisher{Amer. Math. Soc.}, \blocation{Providence, RI}.
\bid{mr={1739680}}
\bptok{imsref}%
\end{bbook}
\endbibitem

\bibitem{Diaconis-Janson-graph-limits}
\begin{barticle}[mr]
\bauthor{\bsnm{Diaconis},~\bfnm{Persi}\binits{P.}} \AND
  \bauthor{\bsnm{Janson},~\bfnm{Svante}\binits{S.}}
(\byear{2008}).
\btitle{Graph limits and exchangeable random graphs}.
\bjournal{Rend. Mat. Appl. (7)}
\bvolume{28}
\bpages{33--61}.
\bid{issn={1120-7183}, mr={2463439}}
\bptok{imsref}%
\end{barticle}
\endbibitem

\bibitem{Easley-Kleinberg-networks-crowds-and-markets}
\begin{bbook}[mr]
\bauthor{\bsnm{Easley},~\bfnm{David}\binits{D.}} \AND
  \bauthor{\bsnm{Kleinberg},~\bfnm{Jon}\binits{J.}}
(\byear{2010}).
\btitle{Networks, Crowds, and Markets: Reasoning About a Highly Connected
  World}.
\bpublisher{Cambridge Univ. Press}, \blocation{Cambridge}.
\bid{mr={2677125}}
\bptok{imsref}%
\end{bbook}
\endbibitem

\bibitem{Faust-Skvoretz-comparing-networks}
\begin{barticle}[author]
\bauthor{\bsnm{Faust},~\bfnm{Katherine}\binits{K.}} \AND
  \bauthor{\bsnm{Skvoretz},~\bfnm{John}\binits{J.}}
(\byear{2002}).
\btitle{Comparing networks across space and time, size and species}.
\bjournal{Sociological Methodology}
\bvolume{32}
\bpages{267--299}.
\bptok{imsref}%
\end{barticle}
\endbibitem

\bibitem{Frank-Strauss86}
\begin{barticle}[mr]
\bauthor{\bsnm{Frank},~\bfnm{Ove}\binits{O.}} \AND
  \bauthor{\bsnm{Strauss},~\bfnm{David}\binits{D.}}
(\byear{1986}).
\btitle{Markov graphs}.
\bjournal{J. Amer. Statist. Assoc.}
\bvolume{81}
\bpages{832--842}.
\bid{issn={0162-1459}, mr={0860518}}
\bptok{imsref}%
\end{barticle}
\endbibitem

\bibitem{Goldenberg-et-al-on-networks}
\begin{barticle}[author]
\bauthor{\bsnm{Goldenberg},~\bfnm{Anna}\binits{A.}},
  \bauthor{\bsnm{Zheng},~\bfnm{Alice~X.}\binits{A.~X.}},
  \bauthor{\bsnm{Fienberg},~\bfnm{Stephen~E.}\binits{S.~E.}} \AND
  \bauthor{\bsnm{Airoldi},~\bfnm{Edoardo~M.}\binits{E.~M.}}
(\byear{2009}).
\btitle{A survey of statistical network models}.
\bjournal{Foundations and Trends in Machine Learning}
\bvolume{2}
\bpages{1--117}.
\bptok{imsref}%
\end{barticle}
\endbibitem

\bibitem{Gondal-ERGM-for-knowledge-production}
\begin{barticle}[author]
\bauthor{\bsnm{Gondal},~\bfnm{Neha}\binits{N.}}
(\byear{2011}).
\btitle{The local and global structure of knowledge production in an emergent
  research field: An exponential random graph analysis}.
\bjournal{Social Networks}
\bvolume{33}
\bpages{20--30}.
\bptok{imsref}%
\end{barticle}
\endbibitem

\bibitem{Gonzalez-Bailon-ergm-for-WWW}
\begin{barticle}[author]
\bauthor{\bsnm{Gonzalez-Bailon},~\bfnm{Sandra}\binits{S.}}
(\byear{2009}).
\btitle{Opening the black box of link formation: Social factors underlying the
  structure of the web}.
\bjournal{Social Networks}
\bvolume{31}
\bpages{271--280}.
\bptok{imsref}%
\end{barticle}
\endbibitem

\bibitem{Goodereau-Kitts-Morris-birds-of-a-feather}
\begin{barticle}[author]
\bauthor{\bsnm{Goodreau},~\bfnm{Steven~M.}\binits{S.~M.}},
  \bauthor{\bsnm{Kitts},~\bfnm{James~A.}\binits{J.~A.}} \AND
  \bauthor{\bsnm{Morris},~\bfnm{Martina}\binits{M.}}
(\byear{2009}).
\btitle{Birds of a feather, or friend of a friend?: Using exponential random
  graph models to investigate adolescent social networks}.
\bjournal{Demography}
\bvolume{46}
\bpages{103--125}.
\bptok{imsref}%
\end{barticle}
\endbibitem

\bibitem{Grunwald-on-MDL}
\begin{bbook}[author]
\bauthor{\bsnm{Gr{\"u}nwald},~\bfnm{Peter~D.}\binits{P.~D.}}
(\byear{2007}).
\btitle{The Minimum Description Length Principle}.
\bpublisher{MIT Press}, \blocation{Cambridge, MA}.
\bptok{imsref}%
\end{bbook}
\endbibitem

\bibitem{Handcock-Gile-sampled-networks}
\begin{barticle}[mr]
\bauthor{\bsnm{Handcock},~\bfnm{Mark~S.}\binits{M.~S.}} \AND
  \bauthor{\bsnm{Gile},~\bfnm{Krista~J.}\binits{K.~J.}}
(\byear{2010}).
\btitle{Modeling social networks from sampled data}.
\bjournal{Ann. Appl. Stat.}
\bvolume{4}
\bpages{5--25}.
\bid{doi={10.1214/08-AOAS221}, issn={1932-6157}, mr={2758082}}
\bptok{imsref}%
\end{barticle}
\endbibitem

\bibitem{statnet-special-issue}
\begin{barticle}[author]
\bauthor{\bsnm{Handcock},~\bfnm{Mark~S.}\binits{M.~S.}},
  \bauthor{\bsnm{Hunter},~\bfnm{David~R.}\binits{D.~R.}},
  \bauthor{\bsnm{Butts},~\bfnm{Carter~T.}\binits{C.~T.}},
  \bauthor{\bsnm{Goodreau},~\bfnm{Steven~M.}\binits{S.~M.}} \AND
  \bauthor{\bsnm{Morris},~\bfnm{Martina}\binits{M.}}
(\byear{2008}).
\btitle{\texttt{statnet}: Software tools for the representation, visualization,
  analysis and simulation of network data}.
\bjournal{Journal of Statistical Software}
\bvolume{24}
\bpages{1--11}.
\bnote{Special issue on \texttt{statnet}}.
\bptok{imsref}%
\end{barticle}
\endbibitem

\bibitem{Hanneke-Fu-Xing-temporal-social-networks}
\begin{barticle}[mr]
\bauthor{\bsnm{Hanneke},~\bfnm{Steve}\binits{S.}},
  \bauthor{\bsnm{Fu},~\bfnm{Wenjie}\binits{W.}} \AND
  \bauthor{\bsnm{Xing},~\bfnm{Eric~P.}\binits{E.~P.}}
(\byear{2010}).
\btitle{Discrete temporal models of social networks}.
\bjournal{Electron. J. Stat.}
\bvolume{4}
\bpages{585--605}.
\bid{doi={10.1214/09-EJS548}, issn={1935-7524}, mr={2660534}}
\bptok{imsref}%
\end{barticle}
\endbibitem

\bibitem{Holland-Leinhardt-p1}
\begin{barticle}[mr]
\bauthor{\bsnm{Holland},~\bfnm{Paul~W.}\binits{P.~W.}} \AND
  \bauthor{\bsnm{Leinhardt},~\bfnm{Samuel}\binits{S.}}
(\byear{1981}).
\btitle{An exponential family of probability distributions for directed
  graphs}.
\bjournal{J. Amer. Statist. Assoc.}
\bvolume{76}
\bpages{33--65}.
\bid{issn={0162-1459}, mr={0608176}}
\bptnote{check related}%
\bptok{imsref}%
\end{barticle}
\endbibitem

\bibitem{Jona-Lasinio-RG-and-prob-theory}
\begin{barticle}[mr]
\bauthor{\bsnm{Jona-Lasinio},~\bfnm{G.}\binits{G.}}
(\byear{2001}).
\btitle{Renormalization group and probability theory}.
\bjournal{Phys. Rep.}
\bvolume{352}
\bpages{439--458}.
\bid{doi={10.1016/S0370-1573(01)00042-4}, issn={0370-1573}, mr={1862625}}
\bptok{imsref}%
\end{barticle}
\endbibitem

\bibitem{Kallenberg-mod-prob}
\begin{bbook}[mr]
\bauthor{\bsnm{Kallenberg},~\bfnm{Olav}\binits{O.}}
(\byear{2002}).
\btitle{Foundations of Modern Probability},
\bedition{2nd} ed.
\bpublisher{Springer}, \blocation{New York}.
\bid{mr={1876169}}
\bptok{imsref}%
\end{bbook}
\endbibitem

\bibitem{Karrer-MEJN-blockmodels-and-community-structure}
\begin{barticle}[mr]
\bauthor{\bsnm{Karrer},~\bfnm{Brian}\binits{B.}} \AND
  \bauthor{\bsnm{Newman},~\bfnm{M.~E.~J.}\binits{M.~E.~J.}}
(\byear{2011}).
\btitle{Stochastic blockmodels and community structure in networks}.
\bjournal{Phys. Rev. E (3)}
\bvolume{83}
\bpages{016107, 10}.
\bid{doi={10.1103/PhysRevE.83.016107}, issn={1539-3755}, mr={2788206}}
\bptok{imsref}%
\end{barticle}
\endbibitem

\bibitem{Kolaczyk-on-network-data}
\begin{bbook}[mr]
\bauthor{\bsnm{Kolaczyk},~\bfnm{Eric~D.}\binits{E.~D.}}
(\byear{2009}).
\btitle{Statistical Analysis of Network Data: Methods and Models}.
\bpublisher{Springer}, \blocation{New York}.
\bid{doi={10.1007/978-0-387-88146-1}, mr={2724362}}
\bptok{imsref}%
\end{bbook}
\endbibitem

\bibitem{Kossinets-effects-of-missing-data}
\begin{barticle}[author]
\bauthor{\bsnm{Kossinets},~\bfnm{Gueorgi}\binits{G.}}
(\byear{2006}).
\btitle{Effects of missing data in social networks}.
\bjournal{Social Networks}
\bvolume{28}
\bpages{247--268}.
\bptok{imsref}%
\end{barticle}
\endbibitem

\bibitem{Krivitsky-Handcock-Morris-adjusting-for-network-size}
\begin{barticle}[mr]
\bauthor{\bsnm{Krivitsky},~\bfnm{Pavel~N.}\binits{P.~N.}},
  \bauthor{\bsnm{Handcock},~\bfnm{Mark~S.}\binits{M.~S.}} \AND
  \bauthor{\bsnm{Morris},~\bfnm{Martina}\binits{M.}}
(\byear{2011}).
\btitle{Adjusting for network size and composition effects in
  exponential-family random graph models}.
\bjournal{Stat. Methodol.}
\bvolume{8}
\bpages{319--339}.
\bid{doi={10.1016/j.stamet.2011.01.005}, issn={1572-3127}, mr={2800354}}
\bptok{imsref}%
\end{barticle}
\endbibitem

\bibitem{Landau-Lifshitz}
\begin{bbook}[author]
\bauthor{\bsnm{Landau},~\bfnm{L.~D.}\binits{L.~D.}} \AND
  \bauthor{\bsnm{Lifshitz},~\bfnm{E.~M.}\binits{E.~M.}}
(\byear{1980}).
\btitle{Statistical Physics}.
\bpublisher{Pergamon Press}, \blocation{Oxford}.
\bptok{imsref}%
\end{bbook}
\endbibitem

\bibitem{Lauritzen-sufficiency-and-prediction}
\begin{barticle}[mr]
\bauthor{\bsnm{Lauritzen},~\bfnm{Steffen~L.}\binits{S.~L.}}
(\byear{1974}).
\btitle{Sufficiency, prediction and extreme models}.
\bjournal{Scand. J. Stat.}
\bvolume{1}
\bpages{128--134}.
\bid{issn={0303-6898}, mr={0378162}}
\bptok{imsref}%
\end{barticle}
\endbibitem

\bibitem{Lauritzen-extremal-families-and-suff-stats}
\begin{bbook}[mr]
\bauthor{\bsnm{Lauritzen},~\bfnm{Steffen~L.}\binits{S.~L.}}
(\byear{1988}).
\btitle{Extremal Families and Systems of Sufficient Statistics}.
\bseries{Lecture Notes in Statistics}
\bvolume{49}.
\bpublisher{Springer}, \blocation{New York}.
\bid{doi={10.1007/978-1-4612-1023-8}, mr={0971253}}
\bptok{imsref}%
\end{bbook}
\endbibitem

\bibitem{Lauritzen-exchangeable-Rasch-matrices}
\begin{barticle}[mr]
\bauthor{\bsnm{Lauritzen},~\bfnm{Steffen~L.}\binits{S.~L.}}
(\byear{2008}).
\btitle{Exchangeable {R}asch matrices}.
\bjournal{Rend. Mat. Appl. (7)}
\bvolume{28}
\bpages{83--95}.
\bid{issn={1120-7183}, mr={2463441}}
\bptok{imsref}%
\end{barticle}
\endbibitem

\bibitem{Lovasz-Szegedy-limits-of-dense-graphs}
\begin{barticle}[mr]
\bauthor{\bsnm{Lov{\'a}sz},~\bfnm{L{\'a}szl{\'o}}\binits{L.}} \AND
  \bauthor{\bsnm{Szegedy},~\bfnm{Bal{\'a}zs}\binits{B.}}
(\byear{2006}).
\btitle{Limits of dense graph sequences}.
\bjournal{J. Combin. Theory Ser. B}
\bvolume{96}
\bpages{933--957}.
\bid{doi={10.1016/j.jctb.2006.05.002}, issn={0095-8956}, mr={2274085}}
\bptok{imsref}%
\end{barticle}
\endbibitem

\bibitem{Lubbers-Snijders-comparison-of-ergms}
\begin{barticle}[author]
\bauthor{\bsnm{Lubbers},~\bfnm{Miranda~J.}\binits{M.~J.}} \AND
  \bauthor{\bsnm{Snijders},~\bfnm{Tom A.~B.}\binits{T.~A.~B.}}
(\byear{2007}).
\btitle{A comparison of various approaches to the exponential random graph
  model: A reanalysis of 102 student networks in school classes}.
\bjournal{Social Networks}
\bvolume{29}
\bpages{489--507}.
\bptok{imsref}%
\end{barticle}
\endbibitem

\bibitem{Mandelbrot-sufficiency-and-estimation-in-thermo}
\begin{barticle}[mr]
\bauthor{\bsnm{Mandelbrot},~\bfnm{Beno{\^{\i}}t}\binits{B.}}
(\byear{1962}).
\btitle{The role of sufficiency and of estimation in thermodynamics}.
\bjournal{Ann. Math. Statist.}
\bvolume{33}
\bpages{1021--1038}.
\bid{issn={0003-4851}, mr={0143592}}
\bptok{imsref}%
\end{barticle}
\endbibitem

\bibitem{Birds-of-a-Feather-review}
\begin{barticle}[author]
\bauthor{\bsnm{McPherson},~\bfnm{Miller}\binits{M.}},
  \bauthor{\bsnm{Smith-Lovin},~\bfnm{Lynn}\binits{L.}} \AND
  \bauthor{\bsnm{Cook},~\bfnm{James~M.}\binits{J.~M.}}
(\byear{2001}).
\btitle{Birds of a feather: Homophily in social networks}.
\bjournal{Annual Review of Sociology}
\bvolume{27}
\bpages{415--444}.
\bptok{imsref}%
\end{barticle}
\endbibitem

\bibitem{Nauenberg-critique-of-q-entropy}
\begin{barticle}[author]
\bauthor{\bsnm{Nauenberg},~\bfnm{Michael}\binits{M.}}
(\byear{2003}).
\btitle{Critique of $q$-entropy for thermal statistics}.
\bjournal{Phys. Rev. E}
\bvolume{67}
\bpages{036114}.
\bptok{imsref}%
\end{barticle}
\endbibitem

\bibitem{MEJN-on-networks}
\begin{bbook}[mr]
\bauthor{\bsnm{Newman},~\bfnm{M.~E.~J.}\binits{M.~E.~J.}}
(\byear{2010}).
\btitle{Networks: An Introduction}.
\bpublisher{Oxford Univ. Press}, \blocation{Oxford}.
\bid{doi={10.1093/acprof:oso/9780199206650.001.0001}, mr={2676073}}
\bptok{imsref}%
\end{bbook}
\endbibitem

\bibitem{Orbanz-projective-limit-techniques}
\begin{bmisc}[author]
\bauthor{\bsnm{Orbanz},~\bfnm{Peter}\binits{P.}}
(\byear{2011}).
\bhowpublished{Projective limit techniques in {Bayesian} nonparametrics.
  Unpublished manuscript.}
\bptok{imsref}%
\end{bmisc}
\endbibitem

\bibitem{Park-MEJN-on-two-star}
\begin{barticle}[author]
\bauthor{\bsnm{Park},~\bfnm{Juyong}\binits{J.}} \AND
  \bauthor{\bsnm{Newman},~\bfnm{Mark E.~J.}\binits{M.~E.~J.}}
(\byear{2004}).
\btitle{Solution of the 2-star model of a network}.
\bjournal{Phys. Rev. E (3)}
\bvolume{70}
\bpages{066146}.
\bptok{imsref}%
\end{barticle}
\endbibitem

\bibitem{Park-MEJN-stat-mech-of-networks}
\begin{barticle}[mr]
\bauthor{\bsnm{Park},~\bfnm{Juyong}\binits{J.}} \AND
  \bauthor{\bsnm{Newman},~\bfnm{M.~E.~J.}\binits{M.~E.~J.}}
(\byear{2004}).
\btitle{Statistical mechanics of networks}.
\bjournal{Phys. Rev. E (3)}
\bvolume{70}
\bpages{066117, 13}.
\bid{doi={10.1103/PhysRevE.70.066117}, issn={1539-3755}, mr={2133807}}
\bptok{imsref}%
\end{barticle}
\endbibitem

\bibitem{Park-MEJN-on-clustered}
\begin{barticle}[author]
\bauthor{\bsnm{Park},~\bfnm{Juyong}\binits{J.}} \AND
  \bauthor{\bsnm{Newman},~\bfnm{Mark E.~J.}\binits{M.~E.~J.}}
(\byear{2006}).
\btitle{Solution for the properties of a clustered network}.
\bjournal{Phys. Rev. E (3)}
\bvolume{72}
\bpages{026136}.
\bptok{imsref}%
\end{barticle}
\endbibitem

\bibitem{Rinaldo-Fienberg-Zhou-network-degeneracy}
\begin{barticle}[mr]
\bauthor{\bsnm{Rinaldo},~\bfnm{Alessandro}\binits{A.}},
  \bauthor{\bsnm{Fienberg},~\bfnm{Stephen~E.}\binits{S.~E.}} \AND
  \bauthor{\bsnm{Zhou},~\bfnm{Yi}\binits{Y.}}
(\byear{2009}).
\btitle{On the geometry of discrete exponential families with application to
  exponential random graph models}.
\bjournal{Electron. J. Stat.}
\bvolume{3}
\bpages{446--484}.
\bid{doi={10.1214/08-EJS350}, issn={1935-7524}, mr={2507456}}
\bptok{imsref}%
\end{barticle}
\endbibitem

\bibitem{Rinaldo-Petrovic-Fienberg}
\begin{bmisc}[author]
\bauthor{\bsnm{Rinaldo},~\bfnm{Alessandro}\binits{A.}},
  \bauthor{\bsnm{Petrovi{\'c}},~\bfnm{S.}\binits{S.}} \AND
  \bauthor{\bsnm{Fienberg},~\bfnm{Stephen~E.}\binits{S.~E.}}
(\byear{2011}).
\bhowpublished{Maximum likelihood estimation in network models. Available at
  arXiv:\arxivurl{1105.6145}.}
\bptok{imsref}%
\end{bmisc}
\endbibitem

\bibitem{Robins-et-al-recent-developments-in-pstar}
\begin{barticle}[author]
\bauthor{\bsnm{Robins},~\bfnm{Garry}\binits{G.}},
  \bauthor{\bsnm{Snijders},~\bfnm{Tom}\binits{T.}},
  \bauthor{\bsnm{Wang},~\bfnm{Peng}\binits{P.}},
  \bauthor{\bsnm{Handcock},~\bfnm{Mark}\binits{M.}} \AND
  \bauthor{\bsnm{Pattison},~\bfnm{Philippa}\binits{P.}}
(\byear{2007}).
\btitle{Recent developments in exponential random graph ($p^*$) models for
  social networks}.
\bjournal{Social Networks}
\bvolume{29}
\bpages{192--215}.
\bptok{imsref}%
\end{barticle}
\endbibitem

\bibitem{Schaefer-youth-co-offending-networks}
\begin{barticle}[author]
\bauthor{\bsnm{Schaefer},~\bfnm{David~R.}\binits{D.~R.}}
(\byear{2012}).
\btitle{Youth co-offending networks: An investigation of social and spatial
  effects}.
\bjournal{Social Networks}
\bvolume{34}
\bpages{141--149}.
\bptok{imsref}%
\end{barticle}
\endbibitem

\bibitem{Schervish-theory-of-stats}
\begin{bbook}[mr]
\bauthor{\bsnm{Schervish},~\bfnm{Mark~J.}\binits{M.~J.}}
(\byear{1995}).
\btitle{Theory of Statistics}.
\bpublisher{Springer}, \blocation{New York}.
\bid{doi={10.1007/978-1-4612-4250-5}, mr={1354146}}
\bptok{imsref}%
\end{bbook}
\endbibitem

\bibitem{supp}
\begin{bmisc}[auto]
\bauthor{\bsnm{Shalizi},~\bfnm{Cosma~Rohilla}\binits{C.~R.}} \AND
  \bauthor{\bsnm{Rinaldo},~\bfnm{Alessandro}\binits{A.}}
(\byear{2013}).
\bhowpublished{Supplement to ``Consistency under sampling of exponential random
  graph models.'' DOI:\doiurl{10.1214/12-AOS1044SUPP}.}
\bptok{imsref}%
\end{bmisc}
\endbibitem

\bibitem{Snijders-on-longitudinal-network-data}
\begin{binproceedings}[author]
\bauthor{\bsnm{Snijders},~\bfnm{Tom A.~B.}\binits{T.~A.~B.}}
(\byear{2005}).
\btitle{Models for longitudinal network data}.
In \bbooktitle{Models and Methods in Social Network Analysis}
(\beditor{\bfnm{Peter~J.}\binits{P.~J.}~\bsnm{Carrington}},
  \beditor{\bfnm{John}\binits{J.}~\bsnm{Scott}} \AND
  \beditor{\bfnm{Stanley}\binits{S.}~\bsnm{Wasserman}}, eds.)
\bpages{215--247}.
\bpublisher{Cambridge Univ. Press}, \blocation{Cambridge}.
\bptok{imsref}%
\end{binproceedings}
\endbibitem

\bibitem{Snijders-Pattison-Robins-Handcock-new-specification}
\begin{barticle}[author]
\bauthor{\bsnm{Snijders},~\bfnm{Tom A.~B.}\binits{T.~A.~B.}},
  \bauthor{\bsnm{Pattison},~\bfnm{Philippa~E.}\binits{P.~E.}},
  \bauthor{\bsnm{Robins},~\bfnm{Garry~L.}\binits{G.~L.}} \AND
  \bauthor{\bsnm{Handcock},~\bfnm{Mark~S.}\binits{M.~S.}}
(\byear{2006}).
\btitle{New specifications for exponential random graph models}.
\bjournal{Sociological Methodology}
\bvolume{36}
\bpages{99--153}.
\bptok{imsref}%
\end{barticle}
\endbibitem

\bibitem{Stumpf-Wiuf-May-subnets-are-not-scale-free}
\begin{barticle}[author]
\bauthor{\bsnm{Stumpf},~\bfnm{Michael P.~H.}\binits{M.~P.~H.}},
  \bauthor{\bsnm{Wiuf},~\bfnm{Carsten}\binits{C.}} \AND
  \bauthor{\bsnm{May},~\bfnm{Robert~M.}\binits{R.~M.}}
(\byear{2005}).
\btitle{Subnets of scale-free networks are not scale-free: Sampling properties
  of networks}.
\bjournal{Proc. Natl. Acad. Sci. USA}
\bvolume{102}
\bpages{4221--4224}.
\bptok{imsref}%
\end{barticle}
\endbibitem

\bibitem{Touchette-large-dev-and-stat-mech}
\begin{barticle}[mr]
\bauthor{\bsnm{Touchette},~\bfnm{Hugo}\binits{H.}}
(\byear{2009}).
\btitle{The large deviation approach to statistical mechanics}.
\bjournal{Phys. Rep.}
\bvolume{478}
\bpages{1--69}.
\bid{doi={10.1016/j.physrep.2009.05.002}, issn={0370-1573}, mr={2560411}}
\bptok{imsref}%
\end{barticle}
\endbibitem

\bibitem{Vermeij-et-al-ergms-of-segregation-networks}
\begin{barticle}[author]
\bauthor{\bsnm{Vermeij},~\bfnm{Lotte}\binits{L.}}, \bauthor{\bparticle{van}
  \bsnm{Duijin},~\bfnm{Marijtje A.~J.}\binits{M.~A.~J.}} \AND
  \bauthor{\bsnm{Baerveldt},~\bfnm{Chris}\binits{C.}}
(\byear{2009}).
\btitle{Ethnic segregation in context: Social discrimination among native
  {Dutch} pupils and their ethnic minority classmates}.
\bjournal{Social Networks}
\bvolume{31}
\bpages{230--239}.
\bptok{imsref}%
\end{barticle}
\endbibitem

\bibitem{Wainwright-Jordan-graphical-models-and-exponential-families}
\begin{barticle}[author]
\bauthor{\bsnm{Wainwright},~\bfnm{Martin~J.}\binits{M.~J.}} \AND
  \bauthor{\bsnm{Jordan},~\bfnm{Michael~I.}\binits{M.~I.}}
(\byear{2008}).
\btitle{Graphical models, exponential families, and variational inference}.
\bjournal{Foundations and Trends in Machine Learning}
\bvolume{1}
\bpages{1--305}.
\bptok{imsref}%
\end{barticle}
\endbibitem

\bibitem{Wasserman-Pattison-1996}
\begin{barticle}[mr]
\bauthor{\bsnm{Wasserman},~\bfnm{Stanley}\binits{S.}} \AND
  \bauthor{\bsnm{Pattison},~\bfnm{Philippa}\binits{P.}}
(\byear{1996}).
\btitle{Logit models and logistic regressions for social networks. {I}. {A}n
  introduction to {M}arkov graphs and {$p$}}.
\bjournal{Psychometrika}
\bvolume{61}
\bpages{401--425}.
\bid{doi={10.1007/BF02294547}, issn={0033-3123}, mr={1424909}}
\bptok{imsref}%
\end{barticle}
\endbibitem

\bibitem{Wasserman-Robins-dependence-graphs-and-pstar}
\begin{bincollection}[author]
\bauthor{\bsnm{Wasserman},~\bfnm{Stanley}\binits{S.}} \AND
  \bauthor{\bsnm{Robins},~\bfnm{Garry}\binits{G.}}
(\byear{2005}).
\btitle{An introduction to random graphs, dependence graphs, and $p^*$}.
In \bbooktitle{Models and Methods in Social Network Analysis}
(\beditor{\bfnm{Peter~J.}\binits{P.~J.}~\bsnm{Carrington}},
  \beditor{\bfnm{John}\binits{J.}~\bsnm{Scott}} \AND
  \beditor{\bfnm{Stanley}\binits{S.}~\bsnm{Wasserman}}, eds.)
\bpages{148--161}.
\bpublisher{Cambridge Univ. Press}, \blocation{Cambridge, England}.
\bptok{imsref}%
\end{bincollection}
\endbibitem

\bibitem{Xiang-Neville-learning-with-one-network}
\begin{binproceedings}[author]
\bauthor{\bsnm{Xiang},~\bfnm{Rongjing}\binits{R.}} \AND
  \bauthor{\bsnm{Neville},~\bfnm{Jennifer}\binits{J.}}
(\byear{2011}).
\btitle{Relational learning with one network: An asymptotic analysis}.
In \bbooktitle{Proceedings of the $14th$ International Conference on Artificial
  Intelligence and Statistics [AISTATS 2011]}
(\beditor{\bfnm{Geoffrey}\binits{G.}~\bsnm{Gordon}},
  \beditor{\bfnm{David}\binits{D.}~\bsnm{Dunson}} \AND
  \beditor{\bfnm{Miroslav}\binits{M.}~\bsnm{Dud{\'i}k}}, eds.).
\bseries{Journal of Machine Learning Research: Workshops and Conference
  Proceedings}
\bvolume{15}
\bpages{779--788}.
\bpublisher{Clarendon Press}, \blocation{Oxford}.
\bptok{imsref}%
\end{binproceedings}\
\endbibitem

\bibitem{Yeomans}
\begin{bbook}[author]
\bauthor{\bsnm{Yeomans},~\bfnm{Julia~M.}\binits{J.~M.}}
(\byear{1992}).
\btitle{Statistical Mechanics of Phase Transitions}.
\bpublisher{Clarendon Press}, \blocation{Oxford}.
\bptok{imsref}%
\end{bbook}
\endbibitem

\end{thebibliography}
\end{document}